\documentclass[10pt]{article}
\usepackage{theorem}
\usepackage{amssymb}
\usepackage{latexsym}

\theorembodyfont{\itshape}
\newtheorem{Theorem}{Theorem}[section]
\newtheorem{Proposition}{Proposition}[section]

\newtheorem{Lemma}{Lemma}[section]
\newtheorem{Claim}{Claim}[section]

\newtheorem{Corollary}{Corollary}[section]

\theorembodyfont{\rmfamily}
\newtheorem{Definition}{Definition}[section]
\newtheorem{Notation}{Notation}[section]
\newtheorem{Remark}{Remark}[section]

\begin{document}
\title{{\bf Remarks on the second sectional geometric genus of quasi-polarized manifolds and their applications} 
\thanks{This research was partially supported by the Grant-in-Aid for Scientific Research (C)
(No.20540045), Japan Society for the Promotion of Science, Japan.}}
\author{YOSHIAKI FUKUMA}
\date{}

\maketitle
\begin{abstract}
In our previous papers, we investigated a lower bound for the second sectional geometric genus $g_{2}(X,L)$ of $n$-dimensional polarized manifolds $(X,L)$ and by using these, we studied the dimension of global sections of $K_{X}+tL$ with $t\geq 2$.
In this paper, we consider the case where $(X,L)$ is a quasi-polarized manifold.
First we will prove $g_{2}(X,L)\geq h^{1}(\mathcal{O}_{X})$ for the following cases:
(a) $n=3$, $\kappa(X)=-\infty$ and $\kappa(K_{X}+L)\geq 0$.
(b) $n\geq 3$ and $\kappa(X)\geq 0$.
Moreover, by using this inequality, we will study $h^{0}(K_{X}+tL)$ for the case where $(X,L)$ is a quasi-polarized $3$-fold.
\end{abstract}

\section{Introduction}
Let $X$ be a smooth projective variety of dimension $n$ defined over the field of complex numbers and let $L$ be a line bundle on $X$.
Then $(X,L)$ is called a {\it quasi-polarized} (resp. {\it polarized}) {\it manifold}
if $L$ is nef and big (resp. ample).
In \cite{Fukuma04-1}, \cite{Fukuma05} and \cite{Fukuma05-2}, we defined the $i$th sectional geometric genus $g_{i}(X,L)$ of $(X,L)$ for any integer $i$ with $0\leq i\leq n$, and we studied some properties of this invariant.
In particular, we proved that $g_{2}(X,L)\geq h^{1}(\mathcal{O}_{X})$ if $(X,L)$ is a polarized manifold with one of the following cases:
\begin{itemize}
\item [\rm (a)] $n=3$, $\kappa(X)=-\infty$ and $\kappa(K_{X}+L)\geq 0$ (see \cite[Theorme 3.3.1 (2)]{Fukuma05-2}).
\item [\rm (b)] $n\geq 3$ and $\kappa(X)\geq 0$ (see \cite[Theorem 2.3.2]{Fukuma05}).
\end{itemize}

Using these results, we also studied the dimension of global sections of $K_{X}+tL$ with $t\geq 2$ for polarized $3$-folds (see \cite{Fukuma06}, \cite{Fukuma07-2} and \cite{Fukuma08-2}).
\par
In this paper, we consider the case where $(X,L)$ is a quasi-polarized manifold.
This generalization is very important.
When we investigate a polarized manifold $(X,L)$, we sometimes need to take a birational morphism $\mu:\widetilde{X}\to X$.
For example, if $(X,L)$ is a polarized variety such that $X$ has singularities, then, by taking a resolution $\mu:\widetilde{X}\to X$, $(\widetilde{X},\mu^{*}(L))$ is not a polarized manifold but a quasi-polarized manifold, and investigation of $(\widetilde{X},\mu^{*}(L))$ makes possible to find some properties of $(X,L)$.
\par
In this paper, first we will study a lower bound for the second sectional geometric genus $g_{2}(X,L)$ of quasi-polarized manifolds $(X,L)$ for the following cases:
\begin{itemize}
\item [\rm (a)] $n=3$, $\kappa(X)=-\infty$ and $\kappa(K_{X}+L)\geq 0$ (Theorem \ref{MT1}).
\item [\rm (b)] $n\geq 3$ and $\kappa(X)\geq 0$ (Theorems \ref{MT3} and \ref{T-T1}).
\end{itemize}

Using these results, we will investigate the dimension of global sections of $K_{X}+tL$ with $t\geq 2$.
Specifically, we get the following results: Let $(X,L)$ be a quasi-polarized $3$-fold.
\begin{itemize}
\item [\rm (a)]
Assume that $\kappa(K_{X}+2L)\geq 0$.
Then $h^{0}(K_{X}+2L)>0$ holds (Theorem \ref{MT4}).
This is an affirmative answer for the case of $n=3$ in \cite[Conjecture NB]{Fukuma98-2}, which can be 
regarded as a generalization of \cite[Conjecture 7.2.7]{BeSo95} proposed by Beltrametti and Sommese.
This also gives a classification of $(X,L)$ with $h^{0}(K_{X}+2L)=0$ (Corollary \ref{T-CR1}).
\item [\rm (b)]
A classification of $(X,L)$ with $h^{0}(K_{X}+2L)=1$ (Theorem \ref{MT5} (a)).
\item [\rm (c)]
A classification of $(X,L)$ with $h^{0}(K_{X}+3L)=0$ (Theorem \ref{MT4-II} (a)).
\item [\rm (d)]
A classification of $(X,L)$ with $h^{0}(K_{X}+3L)=1$ (Theorem \ref{MT5} (b)).
\item [\rm (e)]
$h^{0}(K_{X}+tL)\geq {t-1\choose 3}$ for $t\geq 4$ (Theorem \ref{MT4-II} (b)).
\item [\rm (f)]
A classification of $(X,L)$ with $h^{0}(K_{X}+tL)={t-1\choose 3}$ for some $t\geq 4$ (Theorem \ref{MT5} (c)).
\end{itemize}

After the first version of this paper has been completed, preprints \cite{Horing09} and \cite{Horing09-2} of H\"oring appeared.
We note that Theorem \ref{MT4} is obtained from \cite[1.5 Theorem]{Horing09} and Proposition \ref{T-P1} (ii) below, which is obtained from \cite{Horing09}.
But Theorems \ref{MT4-II} and \ref{MT5} in this paper are not in \cite{Horing09} and these will be useful when we study the dimension of the global sections of adjoint bundles for higher dimensional varieties.
\par
We use standard notation in algebraic geometry.
\par
The author would like to thank Dr. A. H\"oring for sending his preprint \cite{Horing09-2}.

\section{Preliminaries}

\begin{Lemma}\label{L1}
Let $X$ and $C$ be smooth projective varieties with $\dim X=n$ and $\dim C=1$, and let $L$ be a nef and big line bundle on $X$.
Assume that there exists a fiber space $f: X\to C$ such that $h^{0}(K_{F}+L_{F})\neq 0$ 
for a general fiber $F$ of $f$.
Then $f_{*}(K_{X/C}+L)$ is ample.
\end{Lemma}
\noindent{\em Proof.}
First we note that there exists a natural number $m$ such that $(mL)^{n}-n(mL)^{n-1}F>0$.
Then by \cite[(4.1) Lemma]{Demaily96}, there exists a natural number $k$ such that $\mathcal{O}_{X}(k(mL-F))$ has a nontrivial global section.
Hence we have an injective map $\mathcal{O}_{X}(kF)\to \mathcal{O}(kmL)$.
On the other hand, there exists a line bundle $\mathcal{N}$ on $C$ such that
$\mathcal{O}(kF)=f^{*}(\mathcal{N})$.
Hence by \cite[Corollary 1.9]{E-V90} we see that $f_{*}(K_{X/C}+L)$ is ample and we get the assertion. $\Box$

\begin{Lemma}\label{L2}
Let $X$ be a smooth projective variety of dimension $n\geq 2$ and let $V$ be a normal projective variety of dimension $n\geq 2$ such that $V$ has only $\mathbb{Q}$-factorial terminal singularities.
Let $\pi: X\to V$ be a birational morphism such that $X\backslash \pi^{-1}(\mbox{\rm Sing}(V))\cong V\backslash \mbox{\rm Sing}(V)$.
Let $E$ be a $\pi$-exceptional irreducible and reduced divisor on $X$, $A$ a line bundle on $X$ and $L_{1}, \dots , L_{n-2}$ line bundles on $V$.
Then $EA(\pi^{*}(L_{1}))\cdots (\pi^{*}(L_{n-2}))=0$.
\end{Lemma}
\noindent
{\em Proof.}
By \cite[Proposition 4 in section 2, chapter I]{Kleiman66}, we have
$$EA(\pi^{*}(L_{1}))\cdots (\pi^{*}(L_{n-2}))=A|_{E}(\pi^{*}(L_{1}))|_{E}\cdots (\pi^{*}(L_{n-2}))|_{E}.$$
On the other hand since $\dim \mbox{Sing}V\leq n-3$, we have $\dim\pi(E)\leq n-3$.
Here we set $Z:=\pi(E)$.
Then
$$A|_{E}(\pi^{*}(L_{1}))|_{E}\cdots (\pi^{*}(L_{n-2}))|_{E}
=A|_{E}((\pi|_{E})^{*}(L_{1}|_{Z}))\cdots ((\pi|_{E})^{*}(L_{n-2}|_{Z})).$$
Next we consider these intersection numbers.
Here we set $f(t_{1}, \dots , t_{n-1}):=\chi(E, ((\pi|_{E})^{*}(L_{1}|_{Z}))^{\otimes t_{1}}\otimes \cdots \otimes ((\pi|_{E})^{*}(L_{n-2}|_{Z}))^{\otimes t_{n-2}}\otimes (A|_{E})^{\otimes t_{n-1}})$.
Then $f(t_{1}, \dots , t_{n-1})$ is a polynomial of $t_{1}, \dots , t_{n-1}$ of degree at most $n-1$.
Let $C_{1}$ (resp. $C_{2}$) be the coefficient of $t_{1}\cdots t_{n-2}$ (resp. $t_{1}\cdots t_{n-2}t_{n-1}$) in $f(t_{1}, \dots , t_{n-1})$.
Then $f(t_{1}, \dots , t_{n-2}, 0)=\chi(E, ((\pi|_{E})^{*}(L_{1}|_{Z}))^{\otimes t_{1}}\otimes \cdots \otimes ((\pi|_{E})^{*}(L_{n-2}|_{Z}))^{\otimes t_{n-2}})$.
Here we set $g(t_{1}, \dots , t_{n-2}):=f(t_{1}, \dots , t_{n-2}, 0)$.
Then the coefficient of $t_{1}\cdots t_{n-2}$ in $g(t_{1}, \dots , t_{n-2})$ 
is equal to $C_{1}$.
On the other hand since the degree of $g(t_{1}, \dots , t_{n-2})$ is less than $n-2$
(see the proof of \cite[Proposition 6 in section 2, chapter I]{Kleiman66}), we have $C_{1}=0$.
\par
Next we consider $f(t_{1}, \dots , t_{n-2}, 1)$.
Then the coefficient of $t_{1}\cdots t_{n-2}$ in $f(t_{1}, \dots , t_{n-2}, 1)$ is $C_{1}+C_{2}$.
Moreover $f(t_{1}, \dots , t_{n-2}, 1)=\chi(E, ((\pi|_{E})^{*}(L_{1}|_{Z}))^{\otimes t_{1}}\otimes \cdots \otimes ((\pi|_{E})^{*}(L_{n-2}|_{Z}))^{\otimes t_{n-2}}\otimes (A|_{E}))$ and the degree of this polynomial is less than $n-2$ by using the proof of \cite[Proposition 6 in section 2, chapter I]{Kleiman66}.
Hence $C_{1}+C_{2}=0$.
Therefore $C_{2}=0$ since $C_{1}=0$.
Namely the coefficient of $t_{1}\cdots t_{n-1}$ in $f(t_{1},\dots , t_{n-1})$ is zero.
Therefore by the definition of intersection numbers (see \cite{Kleiman66}) we have
$A|_{E}(\pi|_{E})^{*}(L_{1}|_{Z})\cdots (\pi|_{E})^{*}(L_{n-2}|_{Z})=0$.
Hence we get the assertion. $\Box$

\begin{Proposition}\label{P1}
Let $(X,L)$ be a quasi-polarized manifold with $\dim X=n$.
Let $m$ be a positive integer.
Assume that $n\leq 2$ and $\kappa(K_{X}+mL)\geq 0$.
Then $h^{0}(K_{X}+mL)>0$.
\end{Proposition}
\noindent
{\em Proof.}
By the same argument as in the proof of \cite[Theorem 2.8]{Fukuma07-2}, we get the assertion. $\Box$

\begin{Definition}\label{B2}
Let $(X_{1},L_{1})$ and $(X_{2},L_{2})$ be quasi-polarized varieties.
Then $(X_{1},L_{1})$ and $(X_{2},L_{2})$ are said to be {\it birationally equivalent}
if there is another variety $G$ with birational morphisms $g_{i}:G\to X_{i}$
$(i=1,2)$ such that $g_{1}^{*}L_{1}=g_{2}^{*}L_{2}$.
\end{Definition}

\begin{Proposition}\label{T-P1}
Let $(X,L)$ be a quasi-polarized manifold of dimension $n$.
\begin{itemize}
\item [\rm (i)] If $K_{X}+(n-1)L$ is not pseudoeffective, then $(X,L)$ satisfies one of the following.
\begin{itemize}
\item [\rm (i.1)] $g(X,L)=\Delta(X,L)=0$.
Here $g(X,L)$ {\rm (}resp. $\Delta(X,L)${\rm )} denotes the sectional genus {\rm (}resp. the $\Delta$-genus{\rm )} of $(X,L)$.
\item [\rm (i.2)] $(X,L)$ is birationally equivalent to a scroll over a smooth curve.
\end{itemize}
\item [\rm (ii)] If $K_{X}+(n-1)L$ is pseudoeffective, then there exist a quasi-polarized variety $(X^{\prime},L^{\prime})$ which is birationally equivalent to $(X,L)$ such that $X^{\prime}$ is a normal projective variety with only $\mathbb{Q}$-factorial terminal singularities and $K_{X^{\prime}}+(n-1)L^{\prime}$ is nef.
\end{itemize}
\end{Proposition}
\noindent
{\em Proof.}
First we note the following.
\begin{Claim}\label{T-CL1}
$K_{X}+(n-1)L$ is generically nef if and only if $K_{X}+(n-1)$ is pseudoeffective.
\end{Claim}
\noindent
{\em Proof.}
By definition we see that $K_{X}+(n-1)L$ is generically nef 
if $K_{X}+(n-1)L$ is pseudoeffective.
On the other hand by \cite[1.2 Theorem]{Horing09} we can prove that
$K_{X}+(n-1)L$ is pseudoeffective if $K_{X}+(n-1)L$ is generically nef.
Therefore we get the assertion of Claim \ref{T-CL1}. $\Box$
\\
\par
By Claim \ref{T-CL1} we get (i) from \cite[1.3 Proposition]{Horing09-2}.
Moreover we can also prove (ii) by the same argument as Step 1 of Case IV in the proof of
\cite[1.2 Theorem]{Horing09}. $\Box$

\begin{Proposition}\label{T-P2}
Let $(X,L)$ be a quasi-polarized manifold of dimension $n$.
\begin{itemize}
\item [\rm (i)] $g(X,L)\geq 0$ holds.
\item [\rm (ii)] If $g(X,L)=0$, then $\Delta(X,L)=0$.
\item [\rm (iii)] If $g(X,L)=1$, then there exist a quasi-polarized variety 
$(X^{\prime},L^{\prime})$ which is birationally equivalent to $(X,L)$ such that $(X^{\prime},L^{\prime})$ is one of the following two types.
\begin{itemize}
\item [\rm (iii.1)]
$X^{\prime}$ is a normal projective variety with only Gorenstein $\mathbb{Q}$-factorial terminal singularities and $\mathcal{O}(K_{X^{\prime}}+(n-1)L^{\prime})=\mathcal{O}_{X^{\prime}}$ holds.
\item [\rm (iii.2)]
A scroll over a smooth elliptic curve.
\end{itemize}
\end{itemize}
\end{Proposition}
\noindent
{\em Proof.}
For the proof of (i) and (ii), see \cite[Theorems 1.1 and 1.2]{Horing09-2}.
Here we prove (iii).
Assume that $K_{X}+(n-1)L$ is pseudoeffective.
Then by Proposition \ref{T-P1} (ii) we see that there exist a quasi-polarized variety $(X^{\prime},L^{\prime})$, a smooth projective variety $M$ and birational morphisms $\mu_{1}: M\to X$ and $\mu_{2}: M\to X^{\prime}$ such that $X^{\prime}$ is a normal projective variety with only $\mathbb{Q}$-factorial terminal singularities, $\mu_{1}^{*}(L)=\mu_{2}^{*}(L^{\prime})$ and $K_{X^{\prime}}+(n-1)L^{\prime}$ is nef.
Since $g(X^{\prime},L^{\prime})=g(X,L)=1$, we have $(K_{X^{\prime}}+(n-1)L^{\prime})(L^{\prime})^{n-1}=0$.
By the base point free theorem, there exists a natural number $m$ such that $m(K_{X^{\prime}}+(n-1)L^{\prime})$ is free.
So we get $\mathcal{O}(m(K_{X^{\prime}}+(n-1)L^{\prime}))=\mathcal{O}_{X^{\prime}}$.
Namely $K_{X^{\prime}}+(n-1)L^{\prime}$ is numerically trivial.
By the same argument as in the proof of \cite[(3.9) Corollary]{Fujita89-3}, we can prove that
$X^{\prime}$ is Gorenstein and $\mathcal{O}(K_{X^{\prime}}+(n-1)L^{\prime})=\mathcal{O}_{X^{\prime}}$.
\par
Next we consider the case where $K_{X}+(n-1)L$ is not pseudoeffective.
Then by Proposition \ref{T-P1} (i) we see that $(X^{\prime},L^{\prime})$ is a scroll over a smooth curve $N$ because $g(X^{\prime},L^{\prime})=g(X,L)=1$.
In this case, we can easily show that $g(X^{\prime},L^{\prime})=g(N)$.
Hence $N$ is a smooth elliptic curve.
\par
This completes the proof. $\Box$

\begin{Theorem}\label{T-T2}
Let $(f,X,C,L)$ be a quasi-polarized fiber space such that $X$ and $C$ are smooth with $\dim X=n$ and $\dim C=1$.
Then $g(X,L)\geq g(C)$.
\end{Theorem}
\noindent
{\em Proof.}
See \cite[Theorem 3.1]{Fukuma10-2}. $\Box$

\section{Review on the sectional geometric genus and its related topics}

In this section, we will review the $i$th sectional geometric genus of quasi-polarized varieties $(X,L)$ for every integer $i$ with $0\leq i\leq \dim X$.
Up to now, many investigations of $(X,L)$ via the sectional genus were given.
In order to analyze $(X,L)$ more deeply, the author extended these notions.
First in \cite[Definition 2.1]{Fukuma04-1} we gave an invariant called the $i$th sectional geometric genus which can be considered as a generalization of the sectional genus.
Here we recall the definition of this invariant.

\begin{Notation}\label{2-N1}
Let $(X,L)$ be a quasi-polarized variety of dimension $n$, and let $\chi(tL)$ be the Euler-Poincar\'e characteristic of $tL$.
Then $\chi(tL)$ is a polynomial in $t$ of degree $n$, and we set 
$$\chi(tL)=\sum_{j=0}^{n}\chi_{j}(X,L){t+j-1\choose j}.$$
\end{Notation}

\begin{Definition}\label{2-D1}
(\cite[Definition 2.1]{Fukuma04-1} and \cite[Definition 2.1]{Fukuma05-2}.)
Let $(X,L)$ be a quasi-polarized variety of dimension $n$.
\begin{itemize}
\item [\rm (a)] For any integer $i$ with $0\leq i\leq n$ the {\it $i$th sectional H-arithmetic genus $\chi_{i}^{H}(X,L)$ of $(X,L)$} is defined by the following.
$$\chi_{i}^{H}(X,L)=\chi_{n-i}(X,L).$$
\item [\rm (b)] For any integer $i$ with $0\leq i\leq n$ the {\it $i$th sectional geometric genus $g_{i}(X,L)$ of $(X,L)$} is defined by the following.
$$g_{i}(X,L)=(-1)^{i}(\chi_{i}^{H}(X,L)-\chi(\mathcal{O}_{X}))+\sum_{j=0}^{n-i}(-1)^{n-i-j}h^{n-j}(\mathcal{O}_{X}).$$
\end{itemize}
\end{Definition}

\begin{Remark}\label{2-R1} 
\begin{itemize}
\item [\rm (i)] If $i=0$, 
then $\chi_{0}^{H}(X,L)$ and $g_{0}(X,L)$ are equal to the degree $L^{n}$.
If $i=1$, then $g_{1}(X,L)$ is equal to the sectional genus $g(X,L)$ of $(X,L)$.
\item [\rm (ii)] If $i=n$, then $\chi_{n}^{H}(X,L)=\chi(\mathcal{O}_{X})$ and $g_{n}(X,L)=h^{n}({\cal O}_{X})$.
\item [\rm (iii)]
For every integer $i$ with $1\leq i\leq n$ we have
$$\chi_{i}^{H}(X,L)=1-h^{1}(\mathcal{O}_{X})+\cdots +(-1)^{i-1}h^{i-1}(\mathcal{O}_{X})+(-1)^{i}g_{i}(X,L).$$
\item [\rm (iv)]
Using intersection numbers, the second sectional geometric genus can be written as follows:
\begin{eqnarray*}
g_{2}(X,L)&=&-1+h^{1}(\mathcal{O}_{X})+\frac{1}{12}(K_{X}+(n-1)L)(K_{X}+(n-2)L)L^{n-2}\\
&&+\frac{1}{12}c_{2}(X)L^{n-2}+\frac{n-3}{24}(2K_{X}+(n-2)L)L^{n-1}.
\end{eqnarray*}
\end{itemize}
\end{Remark}

In order to explain the geometric meaning of the sectional geometric genus,
we need the following notion.

\begin{Definition}\label{2-D2}
Let $(X,L)$ be a quasi-polarized variety of dimension $n$.
Then we say that $L$ has a {\it $k$-ladder} if there exists a sequence of irreducible and reduced subvarieties $X\supset X_{1}\supset \cdots \supset X_{k}$ such that $X_{i}\in |L_{i-1}|$ for $1\leq i\leq k$, where $X_{0}:=X$, $L_{0}:=L$ and $L_{i}:=L|_{X_{i}}$.
Here we note that if $X$ is smooth and $L$ has no base points, then $L$ has a $k$-ladder for every integer $k$ with $1\leq k\leq n-1$ such that $X_{j}$ is smooth for every $1\leq j\leq k$.\end{Definition}

Then the $i$th sectional geometric genus satisfies the following properties.

\begin{Theorem}\label{2-T1}
{\rm (\cite[Propositions 2.1 and 2.3, and Theorem 2.4]{Fukuma04-2})}
Let $X$ be a projective variety of dimension $n\geq 2$ 
and let $L$ be a nef and big line bundle on $X$.
Assume that $h^{t}(-sL)=0$ for every integers $t$ and $s$ with $0\leq t\leq n-1$ and $1\leq s$, and $|L|$ has an $(n-i)$-ladder for an integer $i$ with $1\leq i\leq n$.
Then the $i$th sectional geometric genus has the following properties.
\begin{itemize}
\item [\rm (i)] $g_{i}(X_{j},L_{j})=g_{i}(X_{j+1},L_{j+1})$ for every integer $j$ with $0\leq j\leq n-i-1$. {\rm (}Here we use the notation in Definition {\rm \ref{2-D2}}.{\rm )}
\item [\rm (ii)]
$g_{i}(X,L)\geq h^{i}(\mathcal{O}_{X})$.
\end{itemize}
\end{Theorem}

In particular, from Theorem \ref{2-T1} (i) and Remark \ref{2-R1} (ii) we see that if $(X,L)$ satisfies the assumption in Theorem \ref{2-T1}, then the $i$th sectional geometric genus is the geometric genus of $i$-dimensional projective variety $X_{n-i}$.
This is the reason why we call this invariant 
the $i$th sectional geometric genus.
From Theorem \ref{2-T1} we see that the $i$th sectional geometric genus is expected to have properties similar to those of the geometric genus of $i$-dimensional projective varieties.
For other results concerning the $i$th sectional geometric genus, for example, see \cite{Fukuma04-1}, \cite{Fukuma04-2}, \cite{Fukuma05} and \cite{Fukuma05-2}.
The following result will be used later.

\begin{Theorem}\label{2-T2}
Let $X$ be a projective variety with $\dim X=n$ and let $L$ be a nef and big line bundle on $X$. 
\\
{\rm (i)} For any integer $i$ with $0\leq i\leq n-1$, we have
$$g_{i}(X,L)=\sum_{j=0}^{n-i-1}(-1)^{n-j}{n-i\choose j}\chi(-(n-i-j)L)
             +\sum_{k=0}^{n-i}(-1)^{n-i-k}h^{n-k}(\mathcal{O}_{X}).$$
\noindent
{\rm (ii)} Assume that $X$ is smooth. 
Then for any integer $i$ with $0\leq i\leq n-1$, we have
$$g_{i}(X,L)=\sum_{j=0}^{n-i-1}(-1)^{j}{n-i\choose j}h^{0}(K_{X}+(n-i-j)L)
             +\sum_{k=0}^{n-i}(-1)^{n-i-k}h^{n-k}(\mathcal{O}_{X}).$$
\end{Theorem}
\noindent{\em Proof.}
(i) By the same argument as in the proof of \cite[Theorem 2.2]{Fukuma04-1}, we obtain
\begin{eqnarray*}
\chi_{n-i}(X,L)
&=&\sum_{j=0}^{n-i}(-1)^{n-i-j}{n-i\choose j}\chi(-(n-i-j)L) \\
&=&\sum_{j=0}^{n-i-1}(-1)^{n-i-j}{n-i\choose j}\chi(-(n-i-j)L)
+\chi(\mathcal{O}_{X}).
\end{eqnarray*}
Hence by Definition \ref{2-D1}, we get the assertion. 
\\
(ii) By using the Serre duality and the Kawamata-Viehweg vanishing theorem, we get the assertion from (i). $\Box$

\begin{Proposition}\label{2-P1}
Let $(X,L)$ be a quasi-polarized manifold of dimension $3$ with $h^{0}(K_{X})=0$.
Then $g_{2}(X,L)\geq h^{2}(\mathcal{O}_{X})\geq 0$ holds.
\end{Proposition}
\noindent
{\em Proof.} By Theorem \ref{2-T2} (ii) we have $g_{2}(X,L)=h^{0}(K_{X}+L)-h^{0}(K_{X})+h^{2}(\mathcal{O}_{X})=h^{0}(K_{X}+L)+h^{2}(\mathcal{O}_{X})\geq h^{2}(\mathcal{O}_{X})\geq 0$. $\Box$
\\
\par
In Section \ref{S4} we need the following lemma.

\begin{Lemma}\label{L3}
Let $X$ be a normal projective variety of dimension $n$ and 
let $\delta: X^{\prime}\to X$ be a resolution of $X$ 
such that $X^{\prime}\backslash \delta^{-1}(\mbox{\rm Sing}(X))\cong X\backslash \mbox{\rm Sing}(X)$.
Let $L$ be a nef and big line bundle on $X$.
Then the following hold.
\begin{itemize}
\item [\rm (i)] If $\dim \mbox{\rm Sing}(X)\leq n-i-1$, then for every integer $k$ with $0\leq k\leq i$ we have $\chi_{k}^{H}(X,L)=\chi_{k}^{H}(X^{\prime},\delta^{*}(L))$.
\item [\rm (ii)] If $\dim \mbox{\rm Sing}(X)=n-2$ and $L$ is ample, then $\chi_{2}^{H}(X^{\prime},\delta^{*}(L))\leq \chi_{2}^{H}(X,L)$ holds.
\end{itemize}
\end{Lemma}
\noindent
{\em Proof.}
Here we put $\mathcal{F}_{q}:=R^{q}\delta_{*}\mathcal{O}_{X^{\prime}}$.
Then $\mathcal{F}_{0}=\mathcal{O}_{X}$ and if $q\geq 1$, then by \cite[(4.2.2) in III]{EGA}
(see also \cite[(1.9) Fact in Chapter 0]{Fujita90})
\begin{eqnarray}
\dim \mbox{Supp}\mathcal{F}_{q}
&\leq& \dim\{ x\in X\ |\ \dim\delta^{-1}(x)\geq q\} \label{L3-1}\\
&\leq& \mbox{min}\{ \dim \mbox{Sing}(X), n-q-1\}. \nonumber
\end{eqnarray}
By the Leray spectral sequence we have
\begin{equation}
\chi(X^{\prime}, (\delta^{*}(L))^{\otimes t})=\sum_{q}(-1)^{q}\chi(X,\mathcal{F}_{q}(L^{\otimes t})). \label{L3-2}
\end{equation}
\noindent
\\
(i) By the assumption that
$\dim \mbox{Sing}(X)\leq n-i-1$, 
we see that $\chi_{l}(X,L)=\chi_{l}(X^{\prime},\delta^{*}(L))$
for every integer $l$ with $n-i\leq l\leq n$.
Therefore by the definition of $\chi_{k}^{H}(X,L)$ we get the first assertion.
\\
\\
(ii) Next we consider the second assertion.
Let $\chi(X,\mathcal{F}_{q}(L^{\otimes t}))=\sum_{j\geq 0}a_{q,j}{t+j-1\choose j}$ for any $q$.
Then by (\ref{L3-1}) the coefficient of ${t+n-3\choose n-2}$ in $\sum_{q}(-1)^{q}\chi(X,\mathcal{F}_{q}(L^{\otimes t}))$ is $a_{0,n-2}-a_{1,n-2}$, and by (\ref{L3-2}) and the definition of the sectional H-arithmetic genus, we have $\chi_{2}^{H}(X^{\prime},\delta^{*}(L))=a_{0,n-2}-a_{1,n-2}$.
Here we note that since $\dim\mbox{Sing}(X)=n-2$, we see that $\chi(X,\mathcal{F}_{q}(L^{\otimes t}))$ is a polynomial of degree at most $n-2$ by (\ref{L3-1}).
Since $L$ is ample, by the Serre vanishing theorem, we have $h^{j}(\mathcal{F}_{1}(L^{\otimes t}))=0$ for every positive integer $j$ and $t\gg 0$.
Hence $a_{1,n-2}\geq 0$.
Since $a_{0,n-2}=\chi_{n-2}(X,L)=\chi_{2}^{H}(X,L)$, we get the second assertion. $\Box$ 
\\
\par
The following are used when we consider the dimension of the global sections of adjoint bundles.

\begin{Definition}\label{2-D3}{\rm (\cite[Definitions 3.1 and 3.2]{Fukuma08-2})}
Let $(X,L)$ be a quasi-polarized manifold of dimension $n$ and let $t$ be a positive integer.
\begin{itemize}
\item [\rm (i)]
Let
\begin{eqnarray*}
F_{0}(t)&:=&h^{0}(K_{X}+tL) \\
F_{i}(t)&:=&F_{i-1}(t+1)-F_{i-1}(t) 
\ \ \mbox{for every integer $i$ with $1\leq i\leq n$.}
\end{eqnarray*}
\item [\rm (ii)]
For every integer $i$ with $0\leq i\leq n$, let 
$$A_{i}(X,L):=F_{n-i}(1).$$ 
We call this $A_{i}(X,L)$ the {\it $i$-th Hilbert coefficient of $(X,L)$}.
\end{itemize}
\end{Definition}

In \cite{Fukuma08-2}, we assumed that $L$ is ample.
But the following results are true for the case where $L$ is nef and big by the same argument as \cite{Fukuma08-2}.

\begin{Remark}\label{2-R2}
(A) (\cite[Remark 3.2 (A)]{Fukuma08-2}) The following hold:
\begin{itemize}
\item [\rm (A.1)] $A_{0}(X,L)=L^{n}$.
\item [\rm (A.2)] $A_{n}(X,L)=h^{0}(K_{X}+L)$.
\end{itemize}
\noindent
(B) (\cite[Proposition 3.2]{Fukuma08-2}) For every integer $i$ with $1\leq i\leq n$ we have
$$
A_{i}(X,L)=g_{i}(X,L)+g_{i-1}(X,L)-h^{i-1}(\mathcal{O}_{X}).
$$
\noindent
(C) (\cite[Remark 3.2 (B)]{Fukuma08-2}) Assume that $\mbox{Bs}|L|=\emptyset$. Here we use notation in Definition \ref{2-D2}.
Then
\begin{eqnarray*}
A_{i}(X,L)
&=&g_{i}(X,L)+g_{i-1}(X,L)-h^{i-1}(\mathcal{O}_{X}) \\
&=&h^{i}(\mathcal{O}_{X_{n-i}})+g_{i-1}(X_{n-i},L_{n-i})-h^{i-1}(\mathcal{O}_{X_{n-i}}) \\
&=&h^{0}(K_{X_{n-i}}+L_{n-i}).
\end{eqnarray*}
\end{Remark}

Then the following result holds.

\begin{Theorem}\label{2-T3}
Let $(X,L)$ be a quasi-polarized manifold of dimension $n$
and let $t$ be a positive integer.
Then the following equality holds.
$$h^{0}(K_{X}+tL)=\sum_{j=0}^{n}{t-1\choose n-j}A_{j}(X,L).$$
\end{Theorem}
\noindent{\em Proof.}
In \cite[Theorem 3.1 and Corollary 3.1]{Fukuma08-2}, we proved this result for the case where $L$ is ample.
But the method still works for the case where $L$ is nef and big. $\Box$
\\
\par
This thereom indicates that it is important to study $A_{i}(X,L)$ when we study the value of 
$h^{0}(K_{X}+tL)$.

\begin{Proposition}\label{2-P2}
Let $(X,L)$ be a quasi-polarized manifold of dimension $n$.
\begin{itemize}
\item [\rm (i)] $A_{1}(X,L)\geq 0$ holds.
\item [\rm (ii)] If $A_{1}(X,L)=0$, then there exist a polarized variety $(V,H)$ and a birational morphism $\pi :X\to V$ such that $(V,H)\cong (\mathbb{P}^{n},\mathcal{O}_{\mathbb{P}^{n}}(1))$ and $L=\pi^{*}(H)$.
\item [\rm (iii)] If $A_{1}(X,L)=1$, then $(X,L)$ satisfies one of the following:
\begin{itemize}
\item [\rm (a)] $(X,L)$ is birationally equivalent to a polarized variety $(Y,B)$ which is one of the following types:
\begin{itemize}
\item [\rm (a.1)]
$Y$ is Gorenstein with $K_{Y}=-(n-1)B$ and $B^{n}=1$.
\item [\rm (a.2)]
A scroll over a smooth elliptic curve with $B^{n}=1$.
\end{itemize}
\item [\rm (b)]
There exist a polarized variety $(V,H)$ and a birational morphism $\pi :X\to V$ such that 
$V$ is a {\rm (}possibly singular{\rm )} quadric hypersurface in $\mathbb{P}^{n+1}$, $H=\mathcal{O}_{V}(1)$ and $L=\pi^{*}(H)$.
\end{itemize}
\end{itemize}
\end{Proposition}
\noindent
{\em Proof.}
(i) Since $A_{1}(X,L)=g_{1}(X,L)+L^{n}-1$, $g_{1}(X,L)\geq 0$ by Proposition \ref{T-P2} (i) and $L^{n}\geq 1$, we get
$A_{1}(X,L)\geq 0$.
\\
(ii) If $A_{1}(X,L)=0$, then by the proof of (i) above we see that $g_{1}(X,L)=0$ and $L^{n}=1$.
By Proposition \ref{T-P2} (ii), we have $\Delta(X,L)=0$.
Hence by \cite[(1.1) Theorem]{Fujita89-3}, we infer that there exist a polarized variety $(V,H)$ and a birational morphism $\pi :X\to V$ such that $H$ is very ample, $L=\pi^{*}(H)$ and $\Delta(V,H)=0$.
By \cite[(5.1)]{Fujita90}, we find that $(V,H)\cong (\mathbb{P}^{n},\mathcal{O}_{\mathbb{P}^{n}}(1))$ because $H^{n}=1$.
Therefore we get the assertion of (ii).
\\
(iii) If $A_{1}(X,L)=0$, then by the proof of (i) above $(X,L)$ satisfies 
one of the following types:
\begin{itemize}
\item [\rm (iii.1)]
$g_{1}(X,L)=1$ and $L^{n}=1$.
\item [\rm (iii.2)]
$g_{1}(X,L)=0$ and $L^{n}=2$.
\end{itemize}

First we consider the case of (iii.1).
Then by Proposition \ref{T-P2} (iii), we see that
$(X,L)$ is birationally equivalent to a polarized variety $(Y,B)$ which is one of the following types:
\begin{itemize}
\item [\rm (iii.1.a)]
$Y$ is Gorenstein with $K_{Y}=-(n-1)B$ and $B^{n}=1$.
\item [\rm (iii.1.b)]
A scroll over a smooth elliptic curve with $B^{n}=1$.
\end{itemize}

Next we consider the case of (iii.2).
Then since $L^{n}=2$, by \cite[(5.1)]{Fujita90} we infer that
there exist a polarized variety $(V,H)$ and a birational morphism $\pi :X\to V$ such that $V$ is a (possibly singular) quadric hypersurface in $\mathbb{P}^{n+1}$ and $H=\mathcal{O}_{V}(1)$.

Therefore we get the assertion. $\Box$
\\
\par
In Theorem \ref{MT6} below, we will study $A_{2}(X,L)$ if $\dim X=3$.

\section{Main results}\label{S4}

\begin{Theorem}\label{MT1}
Let $(X,L)$ be a quasi-polarized $3$-fold.
Assume that $\kappa(X)=-\infty$ and $\kappa(K_{X}+L)\geq 0$.
Then $g_{2}(X,L)\geq h^{1}(\mathcal{O}_{X})$.
\end{Theorem}
\noindent{\em Proof.}
First we note that $h^{3}(\mathcal{O}_{X})=0$ in this case.
If $h^{1}(\mathcal{O}_{X})=0$, then $g_{2}(X,L)=h^{0}(K_{X}+L)+h^{2}(\mathcal{O}_{X})\geq 0=h^{1}(\mathcal{O}_{X})$.
Hence we may assume that $h^{1}(\mathcal{O}_{X})>0$.
Let $\alpha: X\to \mbox{Alb}(X)$ be the Albanese map of $X$.
\par
If $\dim\alpha(X)=2$,
then by the same method as in the proof of \cite[Theorem 3.3.1]{Fukuma05-2}, we get the assertion.
So we may assume that $\dim\alpha(X)=1$.
Then $\alpha(X)$ is a smooth curve and $\alpha: X\to \alpha(X)$ is a fiber space, that is, a surjective morphism with connected fibers.
Set $C:=\alpha(X)$.
Then $g(C)\geq 1$.
Assume $h^{0}(K_{F}+L_{F})=0$ for a general fiber $F$ of $\alpha$.
Then we note that the following holds.
$$h^{0}(K_{F}+L_{F})
=g(F,L_{F})-h^{1}(\mathcal{O}_{F})+h^{2}(\mathcal{O}_{F}).$$
Since $\kappa(X)=-\infty$, we have $\kappa(F)=-\infty$.
Hence $h^{2}(\mathcal{O}_{F})=0$ and we have $g(F,L_{F})=h^{1}(\mathcal{O}_{F})$
because $h^{0}(K_{F}+L_{F})=0$.
By \cite[Theorem 3.1]{Fukuma97-1}, we see that $\kappa(K_{F}+L_{F})=-\infty$.
But this is impossible because we assume that $\kappa(K_{X}+L)\geq 0$.
\par
Therefore $h^{0}(K_{F}+L_{F})\neq 0$ and $\alpha_{*}(K_{X/C}+L)$ is ample by Lemma \ref{L1}.
By the same argument as in the proof of \cite[Theorem 3.3.1]{Fukuma05-2}, we get
$h^{0}(K_{X}+L)>h^{0}(K_{F}+L_{F})(g(C)-1)$.
Therefore
\begin{eqnarray*}
g_{2}(X,L)&=&h^{0}(K_{X}+L)+h^{2}(\mathcal{O}_{X})\\
&>&h^{0}(K_{F}+L_{F})(g(C)-1)\\
&\geq&g(C)-1\\
&=&h^{1}(\mathcal{O}_{X})-1.
\end{eqnarray*}
This completes the proof of Theorem \ref{MT1}. $\Box$

\begin{Theorem}\label{MT2}
Let $V$ be a normal projective variety of dimension $n\geq 3$ such that
$V$ has only $\mathbb{Q}$-factorial terminal singularities.
Let $X$ be a smooth projective variety of dimension $n$ with $\kappa(X)\geq 0$.
Assume that a birational morphism $\pi: X\to V$ satisfies
$X\backslash \pi^{-1}(\mbox{\rm Sing}(V))\cong V\backslash \mbox{\rm Sing}(V)$.
Let $H$ be a nef and big line bundle on $V$ such that $K_{V}+sH$ is nef 
for some positive integer $s$, and let
$H_{1}, \dots , H_{n-2}$ be nef and big line bundles on $V$.
Then
\begin{eqnarray*}
&&c_{2}(X)\pi^{*}(H_{1})\cdots \pi^{*}(H_{n-2}) \\
&&\geq -\frac{s(n-1)}{n}K_{X}\pi^{*}(H)\pi^{*}(H_{1})\cdots \pi^{*}(H_{n-2})
-\frac{s^{2}}{n^{2}}{n\choose 2}(\pi^{*}(H))^{2}\pi^{*}(H_{1})\cdots \pi^{*}(H_{n-2}).
\end{eqnarray*}
\end{Theorem}
\noindent
{\em Proof.}
Let $E:=K_{X}-\pi^{*}(K_{V})$ and $B:=\frac{s}{n}\pi^{*}(H)-\frac{1}{n}E$.
Then $E$ is a $\pi$-exceptional effective $\mathbb{Q}$-divisor on $X$ by assumption.
Let $L$ be an ample line bundle on $X$.
Since by Lemma \ref{L2}
\begin{eqnarray*}
&&Bc_{1}(\Omega_{X}\otimes B)\pi^{*}(H_{1})\cdots \pi^{*}(H_{n-2})\\
&&=\left(\frac{s}{n}\pi^{*}(H)-\frac{1}{n}E\right)(K_{X}+nB)\pi^{*}(H_{1})\cdots \pi^{*}(H_{n-2})\\
&&=\left(\frac{s}{n}\pi^{*}(H)-\frac{1}{n}E\right)(\pi^{*}(K_{V}+sH))\pi^{*}(H_{1})\cdots \pi^{*}(H_{n-2})\\
&&=\left(\frac{s}{n}\pi^{*}(H)\right)(\pi^{*}(K_{V}+sH))\pi^{*}(H_{1})\cdots \pi^{*}(H_{n-2})\\
&&>0,
\end{eqnarray*}
we see that 
$$Bc_{1}(\Omega_{X}\otimes B)(t\pi^{*}(H_{1})+L)\cdots (t\pi^{*}(H_{n-2})+L)>0$$
for any sufficiently large positive integer $t$.
So we fix a positive number $t$ which satisfies this inequality.
Let $H_{i}(t):=t\pi^{*}(H_{i})+L$ for every $i$ with $1\leq i\leq n-2$.
Then $H_{i}(t)$ is ample.
\par
Let
$$0=\mathcal{E}_{0}\subset\mathcal{E}_{1}\subset\cdots\subset\mathcal{E}_{l}=\Omega_{X}$$
be the $(D,H_{1}(t), \dots, H_{n-2}(t))$-semistable filtration of $\Omega_{X}$,
where we put $D:=c_{1}(\Omega_{X})+nB$.
Here we note that $D$ is a nef and $(n-2)$-big $\mathbb{Q}$-Cartier divisor on $X$ 
by assumption.
By \cite[Lemma 2.1]{BlGe71}, there exist a smooth projective variety $Y$ of dimension $n$ 
and a finite surjective morphism $f: Y\to X$ such that $f^{*}(B)$ is a Cartier divisor on $Y$.
We put $A:=f^{*}(B)$ and $\mathcal{G}_{i}:=\mathcal{E}_{i}/\mathcal{E}_{i-1}$ for every integer $i$ with $1\leq i\leq l$.
Let $r_{i}:=\mbox{rank}\ \mathcal{G}_{i}$.
\par
By the same argument as in the proof of \cite[Theorem 2.1]{Fukuma05}, we have
\begin{eqnarray}
&&2c_{2}(f^{*}(\Omega_{X})\otimes A)f^{*}(H_{1}(t))\cdots f^{*}(H_{n-2}(t)) \label{T1-1}\\
&&\geq c_{1}(f^{*}(\Omega_{X})\otimes A)^{2}f^{*}(H_{1}(t))\cdots f^{*}(H_{n-2}(t))\nonumber\\
&&\ \ \ +\sum_{i=1}^{l}\frac{r_{i}-1}{r_{i}}c_{1}(f^{*}(\mathcal{G}_{i})
\otimes A)^{2}f^{*}(H_{1}(t))\cdots f^{*}(H_{n-2}(t)) \nonumber\\
&&\ \ \ -\sum_{i=1}^{l}c_{1}(f^{*}(\mathcal{G}_{i})\otimes A)^{2}
f^{*}(H_{1}(t))\cdots f^{*}(H_{n-2}(t)) \nonumber\\
&&=c_{1}(f^{*}(\Omega_{X})\otimes A)^{2}f^{*}(H_{1}(t))\cdots f^{*}(H_{n-2}(t))\nonumber\\
&&\ \ \ -\sum_{i=1}^{l}\frac{1}{r_{i}}c_{1}(f^{*}(\mathcal{G}_{i})\otimes A)^{2}f^{*}(H_{1}(t))
\cdots f^{*}(H_{n-2}(t)). \nonumber
\end{eqnarray}
(See also \cite[(2.1.10) in Theorem 2.1]{Fukuma05}.)
\par
Here we note that since $c_{1}(\Omega_{X})+nB$ is nef,
we have
\begin{equation}
c_{1}(f^{*}(\Omega_{X})\otimes A)^{2}f^{*}(H_{1}(t))\cdots f^{*}(H_{n-2}(t))>0.\label{T1-2}
\end{equation}

For every integer $i$ with $0\leq i\leq l$ we put
$$\alpha_{i}:=\frac{\delta(f^{*}(\mathcal{G}_{i})\otimes A)
c_{1}(f^{*}(\Omega_{X})\otimes A)f^{*}(H_{1}(t))\cdots f^{*}(H_{n-2}(t))}
{c_{1}(f^{*}(\Omega_{X})\otimes A)^{2}f^{*}(H_{1}(t))\cdots f^{*}(H_{n-2}(t))}.$$
\par
Then  we have the following (see \cite[(2.1.12) and (2.1.13)]{Fukuma05})
\begin{eqnarray}
\sum_{i=1}^{l}r_{i}\alpha_{i}
&=&\sum_{i=1}^{l}\frac{c_{1}(\mathcal{G}_{i}\otimes B)
c_{1}(\Omega_{X}\otimes B)H_{1}(t)\cdots H_{n-2}(t)}
{c_{1}(\Omega_{X}\otimes B)^{2}H_{1}(t)\cdots H_{n-2}(t)}\label{T1-3}\\
&=&\frac{c_{1}(\Omega_{X}\otimes B)
c_{1}(\Omega_{X}\otimes B)H_{1}(t)\cdots H_{n-2}(t)}
{c_{1}(\Omega_{X}\otimes B)^{2}H_{1}(t)\cdots H_{n-2}(t)}\nonumber\\
&=&1\nonumber
\end{eqnarray}
and 
\begin{equation}
\alpha_{1}>\cdots >\alpha_{l}. \label{T1-4}
\end{equation}

By the choice of $t$, we have
\begin{equation}
Bc_{1}(\Omega_{X}\otimes B)(t\pi^{*}(H_{1})+L)\cdots (t\pi^{*}(H_{n-2})+L)>0. \label{T1-5}
\end{equation}

Since $\Omega_{X}$ is generically $(H_{1}(t), \dots , H_{n-2}(t))$-semipositive (\cite[Corollary 6.4]{Miyaoka87}), 
we obtain 
\begin{eqnarray}
&&\delta(f^{*}(\mathcal{G}_{l}))c_{1}(f^{*}(\Omega_{X})\otimes A)f^{*}(H_{1}(t))
\cdots f^{*}(H_{n-2}(t)) \label{T1-6}\\
&&=\delta(f^{*}(\mathcal{G}_{l}))c_{1}(f^{*}(\Omega_{X})\otimes f^{*}(B))f^{*}(H_{1}(t))
\cdots f^{*}(H_{n-2}(t)) \nonumber\\
&&=(\deg f)\delta(\mathcal{G}_{l})(c_{1}(\Omega_{X})+nB)H_{1}(t)\cdots H_{n-2}(t) \nonumber\\
&&\geq 0.\nonumber
\end{eqnarray}

From (\ref{T1-5}) and (\ref{T1-6}), we have
\begin{eqnarray*}
&&\delta(f^{*}(\mathcal{G}_{l})\otimes A)c_{1}(f^{*}(\Omega_{X})\otimes A)f^{*}(H_{1}(t))
\cdots f^{*}(H_{n-2}(t))\\
&&>\delta(f^{*}(\mathcal{G}_{l}))c_{1}(f^{*}(\Omega_{X})\otimes A)f^{*}(H_{1}(t))
\cdots f^{*}(H_{n-2}(t)) \\
&&\geq 0.
\end{eqnarray*}

Hence
\begin{equation}
\alpha_{l}\geq 0. \label{T1-7}
\end{equation}
Here we note that by (\ref{T1-3}), (\ref{T1-4}), and (\ref{T1-7})
\begin{equation}
1\geq \alpha_{1}. \label{T1-8}
\end{equation}
On the other hand by the Hodge index theorem
\begin{eqnarray}
&&c_{1}(f^{*}(\mathcal{G}_{i})\otimes A)^{2}f^{*}(H_{1}(t))\cdots f^{*}(H_{n-2}(t)) \label{T1-9}\\
&&=r_{i}^{2}\delta(f^{*}(\mathcal{G}_{i})\otimes A)^{2}
f^{*}(H_{1}(t))\cdots f^{*}(H_{n-2}(t)) \nonumber\\
&&\leq r_{i}^{2}\alpha_{i}^{2}c_{1}(f^{*}(\Omega_{X})\otimes A)^{2}f^{*}(H_{1}(t))
\cdots f^{*}(H_{n-2}(t)). \nonumber
\end{eqnarray}
Therefore by (\ref{T1-1}), (\ref{T1-2}), (\ref{T1-3}), (\ref{T1-4}), (\ref{T1-7}), 
(\ref{T1-8}) and (\ref{T1-9}) we obtain 
\begin{eqnarray*}
&&2c_{2}(f^{*}(\Omega_{X})\otimes A)f^{*}(H_{1}(t))\cdots f^{*}(H_{n-2}(t)) \\
&&\geq \left(1-\sum_{i=1}^{l}r_{i}\alpha_{i}^{2}\right)
c_{1}(f^{*}(\Omega_{X})\otimes A)^{2}f^{*}(H_{1}(t))\cdots f^{*}(H_{n-2}(t)) \\
&&\geq \left\{1-\left(\sum_{i=1}^{l}r_{i}\alpha_{i}\right)\alpha_{1}\right\}
c_{1}(f^{*}(\Omega_{X})\otimes A)^{2}f^{*}(H_{1}(t))\cdots f^{*}(H_{n-2}(t)) \\
&&=\left(1-\alpha_{1}\right)
c_{1}(f^{*}(\Omega_{X})\otimes A)^{2}f^{*}(H_{1}(t))\cdots f^{*}(H_{n-2}(t)) \\
&&\geq 0.
\end{eqnarray*}

Hence
\begin{equation}
c_{2}(\Omega_{X})H_{1}(t)\cdots H_{n-2}(t)\geq -\left( (n-1)c_{1}(\Omega_{X})B
+{n\choose 2}B^{2}\right) H_{1}(t)\cdots H_{n-2}(t). \label{T1-10}
\end{equation}

Here we note that (\ref{T1-10}) holds for any sufficiently large positive integer $t$.
Hence 

\begin{equation}
c_{2}(\Omega_{X})\pi^{*}(H_{1})\cdots \pi^{*}(H_{n-2})\geq -\left( (n-1)c_{1}(\Omega_{X})B
+{n\choose 2}B^{2}\right)\pi^{*}(H_{1})\cdots \pi^{*}(H_{n-2}). \label{T1-11}
\end{equation}

Since $E$ is an effective divisor with $\dim\pi(E)<n-2$, by Lemma \ref{L2} we have 
$$c_{1}(\Omega_{X})B\pi^{*}(H_{1})\cdots \pi^{*}(H_{n-2})
=c_{1}(\Omega_{X})\left(\frac{s}{n}\pi^{*}(H)\right)\pi^{*}(H_{1})\cdots \pi^{*}(H_{n-2})$$
and
$$B^{2}\pi^{*}(H_{1})\cdots \pi^{*}(H_{n-2})=
\left(\frac{s}{n}\pi^{*}(H)\right)^{2}\pi^{*}(H_{1})\cdots \pi^{*}(H_{n-2}).
$$
Therefore we get the assertion of Theorem \ref{MT2}. $\Box$

\begin{Theorem}\label{MT3}
Let $(X,L)$ be a quasi-polarized manifold of dimension $n\geq 3$ with $\kappa(X)\geq 0$.
Then by Proposition {\rm \ref{T-P1}}, there exist a normal projective variety $V$ of dimension $n$ and a nef and big line bundle $H$ on $V$ such that $V$ has only $\mathbb{Q}$-factorial terminal singularities, $(V,H)$ is birationally equivalent to $(X,L)$ and $K_{V}+(n-1)H$ is nef.
Let $\pi:X^{\prime}\to V$ be a resolution of $V$ such that $X^{\prime}\backslash \pi^{-1}(\mbox{\rm Sing}(V))\cong V\backslash \mbox{\rm Sing}(V)$.
Then the following inequality holds.
\begin{eqnarray*}
g_{2}(X,L)
&\geq& -1+h^{1}(\mathcal{O}_{X})+
\frac{1}{12}\pi^{*}(K_{V})(\pi^{*}(K_{V}+(n-1)H))(\pi^{*}(H))^{n-2}\\
&&\ \ +\frac{n^{2}-3n-1}{12n}\pi^{*}(K_{V}+(n-1)H)(\pi^{*}(H))^{n-1}
+\frac{3n-1}{24n}(\pi^{*}(H))^{n}.
\end{eqnarray*}
\end{Theorem}
\noindent
{\em Proof.}
Then there exist a quasi-polarized manifold $(M,A)$, birational morphisms $\pi_{1}:M\to X$ and $\pi_{2}: M\to V$ such that $A=\pi_{1}^{*}(L)=\pi_{2}^{*}(H)$.
Since $\dim \mbox{Sing}(V)\leq n-3$, we have $\chi_{2}^{H}(X,L)=\chi_{2}^{H}(M,A)=\chi_{2}^{H}(V,H)$ by Lemma \ref{L3}.
By assumption, $V$ has only rational singularities. Hence $h^{j}(\mathcal{O}_{V})=h^{j}(\mathcal{O}_{M})$ for every $j$.
Therefore $g_{2}(X,L)=g_{2}(M,A)=g_{2}(V,H)$ by the definition of the second sectional geometric genus.
Let $\pi:X^{\prime}\to V$ be a resolution of $V$ such that $X^{\prime}\backslash \pi^{-1}(\mbox{Sing}(V))\cong V\backslash \mbox{Sing}(V)$.
Since $g_{2}(X^{\prime},\pi^{*}(H))=g_{2}(V,H)$ by the same argument as above,
from Remark \ref{2-R1} (iv) we have
\begin{eqnarray}
g_{2}(X,L)
&=&g_{2}(V,H) \label{MT3-1}\\
&=&g_{2}(X^{\prime},\pi^{*}(H))\nonumber\\
&=&-1+h^{1}(\mathcal{O}_{X^{\prime}})+\frac{1}{12}(K_{X^{\prime}}+(n-1)\pi^{*}(H))(K_{X^{\prime}}+(n-2)\pi^{*}(H))(\pi^{*}(H))^{n-2}\nonumber\\
&&+\frac{1}{12}c_{2}(X^{\prime})(\pi^{*}(H))^{n-2}+\frac{n-3}{24}(2K_{X^{\prime}}+(n-2)\pi^{*}(H))(\pi^{*}(H))^{n-1}.\nonumber
\end{eqnarray}

Since $K_{V}+(n-1)H$ is nef, by Theorem \ref{MT2}, we have
\begin{eqnarray}
&&c_{2}(X^{\prime})(\pi^{*}(H))^{n-2} \label{MT3-2}\\
&&\geq -(n-1)K_{X^{\prime}}\left(\frac{n-1}{n}\pi^{*}(H)\right)(\pi^{*}(H))^{n-2}
-{n\choose 2}\left(\frac{n-1}{n}\pi^{*}(H)\right)^{2}(\pi^{*}(H))^{n-2}. \nonumber
\end{eqnarray}

Here we note that $K_{X^{\prime}}\pi^{*}(H)^{n-1}=\pi^{*}(K_{V})\pi^{*}(H)^{n-1}$ and
$(K_{X^{\prime}})^{2}\pi^{*}(H)^{n-2}=(\pi^{*}(K_{V}))^{2}\pi^{*}(H)^{n-2}$ hold
by Lemma \ref{L2}.
So we get the assertion by using (\ref{MT3-1}) and (\ref{MT3-2}). $\Box$
\\
\par
In particular, we get the following corollary from Theorem \ref{MT3}.

\begin{Corollary}\label{T-C1}
Let $(X,L)$ be a quasi-polarized $n$-fold with $n\geq 4$ and $\kappa(X)\geq 0$.
Then $g_{2}(X,L)\geq h^{1}(\mathcal{O}_{X})$ holds.
\end{Corollary}

But if $\dim X=3$, then we cannot prove $g_{2}(X,L)\geq h^{1}(\mathcal{O}_{X})$ from the inequality in Theorem \ref{MT3}.
So next we consider the case where $\dim X=3$ and $\kappa(X)\geq 0$.

\begin{Lemma}\label{T-L1}
Let $(X,L)$ be a quasi-polarized $3$-fold.
Assume that $\kappa(K_{X}+L)\geq 0$.
Then there exist a quasi-polarized variety $(X^{+},L^{+})$ of dimension three
such that $X^{+}$ is a normal projective variety with only $\mathbb{Q}$-factorial terminal
singularities, $X^{+}$ is birationally equivalent to $X$, $g_{i}(X,L)=g_{i}(X^{+},L^{+})$ for $i=1,2$ and 
$K_{X^{+}}+L^{+}$ is nef.
\end{Lemma}
\noindent{\em Proof.}
(A) By a result of Fujita \cite[(4.2) Theorem]{Fujita89-3} (see also the proof of \cite[Theorem 4.6]{Horing09}), there exist a normal projective variety $M$ of dimension $3$ with only $\mathbb{Q}$-factorial terminal singularities and a nef and big line bundle $A$ on $M$ such that $(X,L)$ and $(M,A)$ are birationally equivalent and $K_{M}+2A$ is nef.
\\
(B) Assume that there exists an irreducible curve $C$ on $M$ 
such that $(K_{M}+2A)C=0$ and $AC>0$.
Then there exists an extremal ray $R$ on $M$ such that $(K_{M}+2A)R=0$ and $AR>0$.
Let $\rho:M\to M^{\prime}$ be the contraction morphism of $R$.
Assume that $\rho$ is not birational.
Then $\dim M^{\prime}\leq 2$ and there exists a $\mathbb{Q}$-Cartier divisor $B$ on $M^{\prime}$ such that $K_{M}+2A=\rho^{*}(B)$ because $(K_{M}+2A)R=0$.
Hence $K_{M}+A=\rho^{*}(B)-A$.
But this is impossible because $\kappa(K_{M}+A)=\kappa(K_{X}+L)\geq 0$ by assumption.
Hence $\rho$ is birational.
By \cite[Theorem 3.1]{Andreatta95} we see that $\rho$ is a blowing up of a smooth point of $M^{\prime}$.
Let $E$ be its exceptional divisor and $A^{\prime}:=\rho_{*}(A)$.
Then $A^{\prime}$ is a nef and big Cartier divisor on $M^{\prime}$ with $A=\rho^{*}(A^{\prime})-E$ and $K_{M}+2A=\rho^{*}(K_{M^{\prime}}+2A^{\prime})$.
\\
(C) By the same argument as in the proof of \cite[Lemma 4.2.17]{BeSo95} (and by using \cite[Lemma 4.7]{Horing09}), we can prove that each exceptional divisors $E_{i}$ of the contraction morphism of the extremal ray $R_{i}$ as in (B) are disjoint.
\\
(D) By contracting all these extremal rays, we get a normal projective variety $Y$ with only $\mathbb{Q}$-factorial terminal singularities, a nef and big Cartier divisor $H$ on $Y$ and a surjective morphism $\mu: M\to Y$ such that $K_{M}+2A=\mu^{*}(K_{Y}+2H)$ and $K_{M}+A=\mu^{*}(K_{Y}+H)+E_{\mu}$, where $E_{\mu}$ is an effective $\mu$-exceptional divisor.
In particular, we see that $\kappa(K_{M}+A)=\kappa(K_{Y}+H)$.
\\
(E) Next we will prove that $K_{Y}+H$ is nef.
Let $\tau$ be the nef value of $(Y,H)$. Then $\tau$ is rational (see e.g. \cite[Theorem 7.34]{Debarre}).
Assume that $\tau >1$.
\begin{Claim}\label{T-CL1}
There exists an irreducible curve $C$ on $Y$ such that $(K_{Y}+\tau H)C=0$ and $HC>0$.
\end{Claim}
\noindent
{\em Proof.}
Since $K_{Y}+\tau H$ is nef, by the base point free theorem (see \cite[Theorem 3-1-1]{KMM85}), there exists an integer $m\gg 0$ such that $\mbox{Bs}|m(K_{Y}+\tau H)|=\emptyset$.
Let $\Phi: Y\to Z$ be the morphism defined by this linear system.
Then we note that $Z$ is a normal projective variety and $\Phi$ has connected fibers.
\par
Assume that $K_{Y}+\tau H$ is ample.
Let 
$$K_{1}=\{ z\in \overline{NE(Y)}\ |\ \Vert z \Vert =1\}.$$
Then $K_{1}$ is compact.
For any $z\in K_{1}$, we set $f(z):=(K_{Y}+\tau H)z$ and $g(z):=Hz$.
Then $f(z)$ is contiuous and positive on $K_{1}$.
Hence $f(z)$ is bounded from below by a positive rational number $a_{1}$.
Moreover $g(z)$ is also continuous and nonnegative on $K_{1}$.
Hence $g(z)$ is bounded from above by a positive rational number $b_{1}$.
This imples that $(K_{Y}+\tau H)-\frac{a_{1}}{b_{1}}H$ is nonnegative on $K_{1}$.
But this is impossible because $\tau$ is nef value.
Hence $K_{Y}+\tau H$ is not ample and there  exists an irreducible curve $C$ on $Y$ such that $\Phi(C)$ is a point.
In particular $(K_{Y}+\tau H)C=0$ in this case.
\par
Next we assume that $HC=0$ for any irreducible curve $C$ on $Y$ with $(K_{Y}+\tau H)C=0$.
By the construction of $\Phi$, there exists an ample line bundle $G$ on $Z$ such that $p(K_{Y}+\tau H)=\Phi^{*}(G)$ for some positive integer $p$.
Let 
$$K_{2}=\{ z\in \overline{NE(Z)}\ |\ \Vert z \Vert =1\}.$$
Then $K_{2}$ is compact.
For any $z\in K_{2}$, we set $h(z):=Gz$.
Then $h(z)$ is contiuous and positive on $K_{2}$.
Hence $h(z)$ is bounded from below by a positive rational number $a_{2}$.
Let $B$ an irreducible curve on $Y$.
If $\Phi_{*}(B)$ is a point, then $\Phi^{*}(G)B=0$ and $HB=0$ by assumption.
Hence $(\Phi^{*}(G)-\frac{a_{2}}{b_{1}}H)B=0$.
If $\Phi_{*}(B)$ is not a point, then $\Phi_{*}(B)\in \overline{NE(Z)}$.
Hence by the choice of $a_{2}$ and $b_{1}$, we have $(\Phi^{*}(G)-\frac{a_{2}}{b_{1}}H)B\geq 0$.
Hence $\Phi^{*}(G)-\frac{a_{2}}{b_{1}}H=(K_{Y}+\tau H)-\frac{a_{2}}{b_{1}}H$ is nef.
But this is impossible because $\tau$ is nef value.
Therefore we see $HC>0$ for some irreducible curve $C$ on $Y$ with $(K_{Y}+\tau H)C=0$, and we get the assertion of Claim \ref{T-CL1}. $\Box$
\\
\par
We go back to the proof of Lemma \ref{T-L1}.
By Claim \ref{T-CL1} we see that there exists an extremal ray $R$ on $Y$ such that $(K_{Y}+\tau H)R=0$ and $HR>0$.
Let $\psi: Y\to Y^{\prime}$ be the contraction morphism of $R$.
Then $H$ is $\psi$-ample and by \cite[Theorem (5.6)]{AnWi97} (see also \cite[Theorem 2.1]{Mella}) and (D) above, we see that $\psi$ is not birational.
In particular $\dim Y^{\prime}\leq 2$ and there exists a $\mathbb{Q}$-Cartier divisor $B^{\prime}$ on $Y^{\prime}$ such that $K_{Y}+\tau H=\psi^{*}(B^{\prime})$ because $(K_{Y}+\tau H)R=0$.
Hence $K_{Y}+H=\rho^{*}(B^{\prime})-(\tau-1)H$.
But since we assume that $\tau>1$, this is impossible because $\kappa(K_{Y}+H)=\kappa(K_{M}+A)=\kappa(K_{X}+L)\geq 0$ by assumption.
Therefore we get $\tau\leq 1$.
\\
(F) Since $H$ is nef and big, we get $h^{i}(K_{Y}+2H)=h^{i}(K_{Y}+H)=0$ for every $i\geq 1$
by Kawamata-Viehweg vanishing theorem (\cite[Theorem 1-2-5]{KMM85}).
Since $Y$ is Cohen-Macaulay, the Serre duality holds.
Hence $\chi(-2H)=-h^{0}(K_{Y}+2H)$ and $\chi(-H)=-h^{0}(K_{Y}+H)$.
We also note that $h^{i}(\mathcal{O}_{X})=h^{i}(\mathcal{O}_{M})=h^{i}(\mathcal{O}_{Y})$ for every $i\geq 0$.
Therefore by Theorem \ref{2-T2} we can easily see that $g_{1}(X,L)=g_{1}(Y,H)$ and $g_{2}(X,L)=g_{2}(Y,H)$.
\\
(G) By setting $X^{+}:=Y$ and $L^{+}:=H$, we get the assertion. $\Box$

\begin{Theorem}\label{T-T1}
Let $(X,L)$ be a quasi-polarized manifold of dimension three.
Assume that $\kappa(X)\geq 0$.
Then $g_{2}(X,L)\geq h^{1}(\mathcal{O}_{X})$.
\end{Theorem}
\noindent
{\em Proof.}
By Lemma \ref{T-L1} we see that there exist a quasi-polarized variety $(X^{+},L^{+})$ of dimension three
such that $X^{+}$ is a normal variety with only $\mathbb{Q}$-factorial terminal
singularities, $g_{2}(X,L)=g_{2}(X^{+},L^{+})$ and $K_{X^{+}}+L^{+}$ is nef.
Let $\nu: \widetilde{X^{+}}\to X^{+}$ be a resolution of $X^{+}$ such that
$$\widetilde{X^{+}}\backslash\nu^{-1}(\mbox{Sing}(X^{+}))\cong X^{+}\backslash\mbox{Sing}(X^{+}).$$
Here we note that $\dim \mbox{Sing}(X^{+})\leq 0$ and $h^{j}(\mathcal{O}_{X^{+}})=h^{j}(\mathcal{O}_{\widetilde{X^{+}}})$ for $j=0, 1$.
Then by Lemma \ref{L3} (i) and Remark \ref{2-R1} (iii) we have $g_{2}(\widetilde{X^{+}}, \nu^{*}(L^{+}))=g_{2}(X^{+},L^{+})$.
Since $g_{2}(X^{+},L^{+})=g_{2}(X,L)$,
we have $g_{2}(X,L)=g_{2}(\widetilde{X^{+}}, \nu^{*}(L^{+}))$.
Here we use Theorem \ref{MT2}.
Then the following inequality holds.
\begin{eqnarray}
c_{2}(\widetilde{X^{+}})(\nu^{*}(L^{+})) 
&\geq& -\frac{2}{3}K_{\widetilde{X}}(\nu^{*}(L^{+}))^{2}-\frac{1}{9}{3\choose 2}(\nu^{*}(L^{+}))^{3} \label{T-T1-a} \\
&=&-\frac{2}{3}K_{\widetilde{X}}(\nu^{*}(L^{+}))^{2}
-\frac{1}{3}(\nu^{*}(L^{+}))^{3}. \nonumber
\end{eqnarray}
Therefore by (\ref{T-T1-a}), Remark \ref{2-R1} (iv) and Lemma \ref{L2} we have
\begin{eqnarray*}
&&g_{2}(\widetilde{X^{+}}, \nu^{*}(L^{+})) \\
&&=-1+h^{1}(\mathcal{O}_{\widetilde{X^{+}}})
+\frac{1}{12}(K_{\widetilde{X^{+}}}+2\nu^{*}(L^{+}))(K_{\widetilde{X^{+}}}+\nu^{*}(L^{+}))(\nu^{*}(L^{+}))\\
&&+\frac{1}{12}c_{2}(\widetilde{X^{+}})(\nu^{*}(L^{+}))\\
&&\geq -1+h^{1}(\mathcal{O}_{\widetilde{X^{+}}})
+\frac{1}{12}((K_{\widetilde{X^{+}}})^{2}+3K_{\widetilde{X^{+}}}\nu^{*}(L^{+})+2\nu^{*}(L^{+})^{2})(\nu^{*}(L^{+}))\\
&&-\frac{1}{18}K_{\widetilde{X^{+}}}\nu^{*}(L^{+})-\frac{1}{36}(\nu^{*}(L^{+}))^{3}\\
&&=-1+h^{1}(\mathcal{O}_{\widetilde{X^{+}}})
+\frac{1}{12}((K_{X^{+}})^{2}+3K_{X^{+}}(L^{+})+2(L^{+})^{2})L^{+}\\
&&-\frac{1}{18}K_{X^{+}}(L^{+})^{2}-\frac{1}{36}(L^{+})^{3}\\
&&=-1+h^{1}(\mathcal{O}_{\widetilde{X^{+}}})
+\frac{1}{12}(K_{X^{+}})^{2}L^{+}+\frac{7}{36}K_{X^{+}}(L^{+})^{2}+\frac{5}{36}(L^{+})^{3}\\
&&=-1+h^{1}(\mathcal{O}_{\widetilde{X^{+}}})
+\frac{1}{12}(K_{X^{+}}+L^{+})K_{X^{+}}L^{+}+\frac{1}{9}(K_{X^{+}}+L^{+})(L^{+})^{2}+\frac{1}{36}(L^{+})^{3}\\
&&>h^{1}(\mathcal{O}_{\widetilde{X^{+}}})-1.
\end{eqnarray*}

So we have
$$
g_{2}(X,L)=g_{2}(\widetilde{X^{+}}, \nu^{*}(L^{+}))
\geq h^{1}(\mathcal{O}_{\widetilde{X^{+}}})
=h^{1}(\mathcal{O}_{X})
$$
and we get the assertion. $\Box$

\begin{Theorem}\label{MT6}
Let $(X,L)$ be a quasi-polarized manifold of dimension $3$.
Then the following hold.
\begin{itemize}
\item [\rm (i)] If $\kappa(X)\geq 0$, then $A_{2}(X,L)\geq 2$.
\item [\rm (ii)] Assume that $\kappa(X)=-\infty$.
\begin{itemize}
\item [\rm (ii.1)] If $h^{1}(\mathcal{O}_{X})=0$, then $A_{2}(X,L)=g_{2}(X,L)+g_{1}(X,L)\geq 0$.
\item [\rm (ii.2)] If $h^{1}(\mathcal{O}_{X})>0$ and the dimension of the image of the Albanese map of $X$ is one, then $A_{2}(X,L)\geq g_{2}(X,L)\geq 0$.
\item [\rm (ii.3)] If $h^{1}(\mathcal{O}_{X})>0$ and the dimension of the image of the Albanese map of $X$ is two, then $A_{2}(X,L)\geq g_{1}(X,L)-1+\chi(\mathcal{O}_{S})\geq 0$, where $S$ is a resolution of the image of the Albanese map of $X$.
\end{itemize}
\end{itemize}
\end{Theorem}
\noindent
{\em Proof.}
(i) By Theorem \ref{T-T1}, we have
$g_{2}(X,L)\geq h^{1}(\mathcal{O}_{X})$.
On the other hand, since $\kappa(X)\geq 0$, we see that
$g_{1}(X,L)=1+(1/2)(K_{X}+2L)L^{2}\geq 2$.
Hence
$A_{2}(X,L)=g_{2}(X,L)+g_{1}(X,L)-h^{1}(\mathcal{O}_{X})\geq 2$.
\\
(ii) Assume that $\kappa(X)=-\infty$.
\\
(ii.1) The case of $h^{1}(\mathcal{O}_{X})=0$.
Then $A_{2}(X,L)=g_{2}(X,L)+g_{1}(X,L)$ by Remark \ref{2-R2} (B).
Since $g_{2}(X,L)\geq 0$ by Proposition \ref{2-P1} and $g_{1}(X,L)\geq 0$ by \cite[(4.8) Corollary]{Fujita89-3} or Proposition \ref{T-P2} (i),
we have $A_{2}(X,L)\geq g_{1}(X,L)\geq 0$.
\\
(ii.2) The case where the dimension of the image of the Albanese map of $X$ is one.
Let $\alpha: X\to \mbox{Alb}(X)$ be the Albanese map of $X$.
Then $\alpha(X)$ is a smooth curve and $\alpha: X\to \alpha(X)$ is a surjective morphism with connected fibers.
Let $C:=\alpha(X)$.
Then by Theorem \ref{T-T2} we have $g_{1}(X,L)\geq g(C)$.
Since $g(C)=h^{1}(\mathcal{O}_{X})$, we get $g_{1}(X,L)\geq h^{1}(\mathcal{O}_{X})$.
Hence $A_{2}(X,L)=g_{2}(X,L)+g_{1}(X,L)-h^{1}(\mathcal{O}_{X})\geq g_{2}(X,L)\geq 0$.
\\
(ii.3) The case where the dimension of the image of the Albanese map of $X$ is two.
Then there exist a smooth projective $3$-fold $X^{\prime}$, a smooth projective surface $S$, 
birational maps $\mu: X^{\prime}\to X$ and $\nu:S\to \alpha(X)$ and a surjective morphism $f:X^{\prime}\to S$ such that $\alpha\circ\mu=\nu\circ f$.
Then we note that 
$h^{1}(\mathcal{O}_{S})=h^{1}(\mathcal{O}_{X})=h^{1}(\mathcal{O}_{X^{\prime}})$ and 
$h^{2}(\mathcal{O}_{X^{\prime}})\geq h^{2}(\mathcal{O}_{S})$ hold.
Therefore 
$$1-h^{1}(\mathcal{O}_{X^{\prime}})+h^{2}(\mathcal{O}_{X^{\prime}}) 
\geq 1-h^{1}(\mathcal{O}_{S})+h^{2}(\mathcal{O}_{S})=\chi(\mathcal{O}_{S}).$$
On the other hand $g_{2}(X,L)\geq h^{2}(\mathcal{O}_{X})$ holds by Proposition \ref{2-P1}.
Hence 
\begin{eqnarray*}
g_{2}(X,L)
&\geq& h^{2}(\mathcal{O}_{X}) \\
&=&h^{2}(\mathcal{O}_{X^{\prime}}) \\
&\geq& h^{1}(\mathcal{O}_{X^{\prime}})-1+\chi(\mathcal{O}_{S}) \\
&=&h^{1}(\mathcal{O}_{X})-1+\chi(\mathcal{O}_{S}).
\end{eqnarray*}
Therefore we have
\begin{eqnarray}
A_{2}(X,L)
&=&g_{2}(X,L)+g_{1}(X,L)-h^{1}(\mathcal{O}_{X}) \label{MT5-0}\\
&\geq&g_{1}(X,L)-1+\chi(\mathcal{O}_{S}).\nonumber
\end{eqnarray}
Here we note that $\chi(\mathcal{O}_{S})\geq 0$ since $\kappa(S)\geq 0$.
We also note that $g_{1}(X,L)\geq 1$ because $h^{1}(\mathcal{O}_{X})=0$ holds if $g_{1}(X,L)=0$ by \cite[(4.8) Corollary and (1.1) Theorem]{Fujita89-3} or Proposition \ref{T-P2} (ii).
Therefore we get the assertion of (ii.3) and these complete the proof of Theorem \ref{MT6}. $\Box$

\begin{Remark}\label{T-R2}
For the case of $\dim X=3$,
we can also prove Theorem \ref{MT6} (i) by using the inequality in Theorem \ref{MT3}.
Here we use notation in Theorem \ref{MT3}.
Then by Theorem \ref{MT3}, we have
\begin{eqnarray*}
A_{2}(X,L)
&=&g_{2}(X,L)+g_{1}(X,L)-h^{1}(\mathcal{O}_{X}) \\
&\geq& -1+\frac{1}{12}\pi^{*}(K_{V})\pi^{*}(K_{V}+2H)\pi^{*}(H)\\
&&+\frac{1}{18}\pi^{*}(H)^{3}-\frac{1}{36}\pi^{*}(K_{V})\pi^{*}(H)^{2}\\
&&+1+\frac{1}{2}\pi^{*}(K_{V}+2H)\pi^{*}(H)^{2} \\
&=&\frac{1}{12}K_{V}(K_{V}+2H)H+\frac{17}{36}K_{V}H^{2}+\frac{19}{18}H^{3}\\
&>&1.
\end{eqnarray*}
(Here we note that $g_{1}(X,L)=1+\frac{1}{2}(K_{M}+2\pi^{*}H)(\pi^{*}(H))^{2}=1+\frac{1}{2}\pi^{*}(K_{V}+2H)(\pi^{*}(H))^{2}$.)
\end{Remark}

The following theorem shows that \cite[Conjecture NB]{Fukuma98-2} for the case of $\dim X=3$ is true, which is a quasi-polarized manifolds' version of a conjecture of Beltrametti and Sommese \cite[Conjecture 7.2.7]{BeSo95}. 

\begin{Theorem}\label{MT4}
Let $(X,L)$ be a quasi-polarized manifold of dimension $3$.
Assume that $\kappa(K_{X}+2L)\geq 0$.
Then $h^{0}(K_{X}+2L)>0$ holds.
\end{Theorem}
\noindent
{\em Proof.}
First we note that $A_{3}(X,L)\geq 0$ and $A_{0}(X,L)\geq 1$ hold in general (see Remark \ref{2-R2} (A)).
Moreover we have $A_{1}(X,L)\geq 0$ because $g_{1}(X,L)\geq 0$ (see \cite[(4.8) Corollary]{Fujita89-3} or Proposition \ref{T-P2} (i)). 
\\
(I) If $\kappa(X)\geq 0$, then by Theorem \ref{MT6} (i) we have $A_{2}(X,L)\geq 2$.
Therefore by Theorem \ref{2-T3} we have $h^{0}(K_{X}+2L)\geq 2$.
\\
(II) Next we assume that $\kappa(X)=-\infty$.
\\
(II.1) If $h^{1}(\mathcal{O}_{X})>0$,
then we take the Albanese map $\alpha: X\to \mbox{Alb}(X)$.
By taking its Stein factorization, if necessary, we make a fiber space $\alpha: X\to Y$ over a normal projective variety $Y$.
Let $F$ be a general fiber of $\alpha$.
Then $\dim F\leq 2$, and $\kappa(K_{F}+2L_{F})\geq 0$ since $\kappa(K_{X}+2L)\geq 0$.
Hence by Proposition \ref{P1}, we have $h^{0}(K_{F}+2L_{F})>0$.
Therefore we have $h^{0}(K_{X}+2L)>0$ by \cite[Lemma 4.1]{ChHa02}.
\\
(II.2) Next we consider the case of $h^{1}(\mathcal{O}_{X})=0$.
Then $h^{0}(K_{X}+2L)=A_{2}(X,L)+A_{3}(X,L)\geq A_{2}(X,L)$.
On the other hand, since $\kappa(K_{X}+2L)\geq 0$, we have
$$g_{1}(X,L)=1+\frac{1}{2}(K_{X}+2L)L^{2}\geq 1.$$
Hence by Proposition \ref{2-P1}
$$A_{2}(X,L)=g_{2}(X,L)+g_{1}(X,L)\geq 1.$$
Therefore we get $h^{0}(K_{X}+2L)\geq 1$.
\par
This completes the proof of Theorem \ref{MT4}. $\Box$

\begin{Remark}
This result is also obtained from \cite[1.5 Theorem]{Horing09} 
and Proposition \ref{T-P1} (ii).
\end{Remark}

By Theorem \ref{MT4} and \cite[(4.2) Theorem]{Fujita89-3}, we get the following result.

\begin{Corollary}\label{T-CR1}
Let $(X,L)$ be a quasi-polarized manifold of dimension $3$.
Then $h^{0}(K_{X}+2L)=0$ if and only if $(X,L)$ is birationally equivalent to a scroll over a smooth curve or a quasi-polarized variety $(V,H)$ such that $V$ is a normal projective variety with only $\mathbb{Q}$-factorial terminal singularities and $\Delta(V,H)=0$.
\end{Corollary}

\begin{Theorem}\label{MT4-II}
Let $(X,L)$ be a quasi-polarized manifold of dimension $3$.
\begin{itemize}
\item [\rm (a)] $h^{0}(K_{X}+3L)=0$ if and only if there exists a birational morphism $f:X\to \mathbb{P}^{3}$ such that $L=f^{*}(\mathcal{O}_{\mathbb{P}^{3}}(1))$. In particular, if $\kappa(K_{X}+3L)\geq 0$, then $h^{0}(K_{X}+3L)\geq 1$.
\item [\rm (b)] For $t\geq 4$, then $h^{0}(K_{X}+tL)\geq {t-1\choose 3}$.
\end{itemize}
\end{Theorem}
\noindent
{\em Proof.}
(a) First we consider $h^{0}(K_{X}+3L)$.
Then by Theorem \ref{2-T3} we have 
\begin{eqnarray*}
h^{0}(K_{X}+3L)
&=&\sum_{k=0}^{3}{2\choose 3-k}A_{k}(X,L) \\
&=&A_{1}(X,L)+2A_{2}(X,L)+A_{3}(X,L).
\end{eqnarray*}
Assume that $\kappa(X)\geq 0$.
Then by Theorem \ref{MT6} we have $A_{2}(X,L)\geq 2$.
We also note that $A_{3}(X,L)\geq 0$ and $A_{1}(X,L)=g_{1}(X,L)+L^{3}-1\geq 2$ because $g_{1}(X,L)=1+(1/2)(K_{X}+2L)L^{2}\geq 2$.
Hence $h^{0}(K_{X}+3L)\geq 6$.
So we may assume that $\kappa(X)=-\infty$.
Here we note that $A_{i}(X,L)\geq 0$ for $i=1,2,3$ 
by Remark \ref{2-R2} (A), Proposition \ref{2-P2} and Theorem \ref{MT6} (ii).
Hence if $h^{0}(K_{X}+3L)=0$, then $A_{3}(X,L)=0$, $A_{2}(X,L)=0$ and $A_{1}(X,L)=0$.
In particular, $A_{1}(X,L)=0$ implies $g_{1}(X,L)=0$ and $L^{3}=1$.
Therefore by \cite[(4.8) Corollary]{Fujita89-3} or Proposition \ref{T-P2} (ii) we see that $\Delta(X,L)=0$.
Moreover by \cite[(1.1) Theorem]{Fujita89-3} we see that there exist a projective variety $W$, a birational morphism $f: X\to W$ and a very ample line bundle $H$ on $W$ such that $L=f^{*}(H)$ and $\Delta(W,H)=0$.
Since $L^{3}=1$, we have $H^{3}=1$.
Therefore $\Delta(W,H)=0$ implies that $h^{0}(H)=4$.
Since $H$ is very ample and $\dim W=3$, we have $W$ is isomorphic to $\mathbb{P}^{3}$.
We can easily check that $h^{0}(K_{X}+3L)=0$ if there exists a birational morphism $f:X\to \mathbb{P}^{3}$ such that $L=f^{*}(\mathcal{O}_{\mathbb{P}^{3}}(1))$.
\\
\\
(b) If $t\geq 4$, then by Theorem \ref{2-T3} we have 
$$h^{0}(K_{X}+tL)=\sum_{k=0}^{3}{t-1\choose 3-k}A_{k}(X,L).$$

By Proposition \ref{2-P2}, Theorem \ref{MT6} and Remark \ref{2-R2} (A) we have
$A_{i}(X,L)\geq 0$ for $i=1,2,3$ and $A_{0}(X,L)\geq 1$.
Therefore we get the assertion. $\Box$

\begin{Theorem}\label{MT5}
Let $(X,L)$ be a quasi-polarized manifold of dimension $3$.
\begin {itemize}
\item [\rm (a)] Assume that $h^{0}(K_{X}+2L)=1$ holds.
Then $(X,L)$ satisfies one of the following three types.
\begin{itemize}
\item [\rm (a.1)] $(X,L)$ is birationally equivalent to $(V,H)$, where $V$ is a normal projective variety with only Gorenstein $\mathbb{Q}$-factorial terminal singularities, $\mathcal{O}(K_{V}+2H)\cong\mathcal{O}_{V}$ and $\Delta(V,H)=1$.
\item [\rm (a.2)] There exist an Abelian surface $S^{\prime}$ and a surjective morphism with connected fibers $f^{\prime}:X\to S^{\prime}$
such that a general fiber $F^{\prime}$ of $f^{\prime}$ is isomorphic to $\mathbb{P}^{1}$ and $L_{F^{\prime}}=\mathcal{O}_{\mathbb{P}^{1}}(1)$.
\item [\rm (a.3)] There exist a smooth elliptic curve $C$ and a surjective morphism with connected fibers $f:X\to C$ such that $L_{F}$-minimalization of $(F,L_{F})$ 
 {\rm (}for the definition of the $L_{F}$-minimalization see {\rm \cite[Definition 1.9 (2)]{Fukuma97-1}}{\rm )}
is isomorphic to either $(\mathbb{P}^{2},\mathcal{O}_{\mathbb{P}^{2}}(2))$ or a scroll over a smooth curve.
\end{itemize}
\item [\rm (b)] Assume that $h^{0}(K_{X}+3L)=1$ holds.
Then $(X,L)$ satisfies one of the following types.
\begin{itemize}
\item [\rm (b.1)]
$(X,L)$ is birationally equivalent to a quasi-polarized variety $(V,H)$ such that $(V,H)$ is one of the following types:
\begin{itemize}
\item [\rm (b.1.1)]
$V$ is a normal projective variety with only $\mathbb{Q}$-factorial terminal singularities, $\omega_{V}\otimes \mathcal{O}(H)^{\otimes 2}\cong\mathcal{O}_{V}$, $H^{3}=1$ and $\Delta(V,H)=1$.
\item [\rm (b.1.2)] $(V,H)$ is a scroll over a smooth elliptic curve and $H^{3}=1$.
\end{itemize}
\item [\rm (b.2)]
There exist a normal projective variety $W$, a very ample line bundle $H$ and a birational morphism $\mu:X \to W$ such that $L=\mu^{*}(H)$, $\Delta(W,H)=0$ and $H^{3}=2$.
\end{itemize}
\item [\rm (c)] Assume that $h^{0}(K_{X}+tL)={t-1\choose 3}$ holds for some $t\geq 4$.
Then there exists a birational morphism $f:X\to \mathbb{P}^{3}$ such that $L=f^{*}(\mathcal{O}_{\mathbb{P}^{3}}(1))$.
\end{itemize}
\end{Theorem}
\noindent
{\em Proof.}
(a) First we assume that $h^{0}(K_{X}+2L)=1$.
If $\kappa(X)\geq 0$, then we have $h^{0}(K_{X}+2L)\geq 2$ 
by the proof of Theorem \ref{MT4}.
So we may assume that $\kappa(X)=-\infty$.
\par
If $\kappa(K_{X}+L)\geq 0$, then by Theorem \ref{MT1} we have $g_{2}(X,L)\geq h^{1}(\mathcal{O}_{X})$.
Hence we see that $A_{2}(X,L)=g_{2}(X,L)+g_{1}(X,L)-h^{1}(\mathcal{O}_{X})\geq g_{1}(X,L)$.
On the other hand, $\kappa(K_{X}+L)\geq 0$ implies that $g_{1}(X,L)\geq 2$.
Hence we have $A_{2}(X,L)\geq 2$ and by Theorem \ref{2-T3} and Remark \ref{2-R2} (A) we get $h^{0}(K_{X}+2L)\geq 2$.
Therefore we see that 
\begin{equation}
\kappa(K_{X}+L)=-\infty. \label{MT5+2}
\end{equation}
In particular $h^{0}(K_{X}+L)=0$.
On the other hand by Remark \ref{2-R2} (A.2) we have $A_{3}(X,L)=h^{0}(K_{X}+L)$.
Hence
\begin{eqnarray}
1&=&h^{0}(K_{X}+2L) \label{MT5+1}\\
&=&A_{2}(X,L)+A_{3}(X,L) \nonumber\\
&=&A_{2}(X,L). \nonumber
\end{eqnarray}
\noindent
\\
(a.1) The case of $h^{1}(\mathcal{O}_{X})=0$.
Then $A_{2}(X,L)\geq 0$ by Theorem \ref{MT6} (ii.1).
On the other hand, since $h^{0}(K_{X}+2L)=A_{2}(X,L)+A_{3}(X,L)$ and $A_{3}(X,L)=h^{0}(K_{X}+L)=0$, we have $g_{1}(X,L)\leq 1$.
Since $\kappa(K_{X}+2L)\geq 0$ by assumption, we have $g_{1}(X,L)=1$.
Since $h^{1}(\mathcal{O}_{X})=0$ in this case, by \cite[(4.9) Corollary]{Fujita89-3} or Proposition \ref{T-P2} (iii) we see that
there exists a quasi-polarized variety $(V,H)$ such that $V$ is a normal projective variety with only Gorenstein $\mathbb{Q}$-factorial terminal singularities, $\mathcal{O}(K_{V}+2H)=\mathcal{O}_{V}$, $\Delta(H)=1$ and $(V,H)$ is birationally equivalent to $(X,L)$.
\\
(a.2) The case of $h^{1}(\mathcal{O}_{X})>0$.
Let $\alpha: X\to \mbox{Alb}(X)$ be the Albanese map of $X$.
\\
(a.2.1) The case of $\dim\alpha(X)=2$.
Then there exist a smooth projective $3$-fold $X^{\prime}$, a smooth projective surface $S$, 
birational maps $\mu: X^{\prime}\to X$ and $\nu:S\to \alpha(X)$ and a surjective morphism $f:X^{\prime}\to S$ such that $\alpha\circ\mu=\nu\circ f$.
Then by Theorem \ref{MT6} (ii.3) we have $A_{2}(X,L)\geq g_{1}(X,L)-1+\chi(\mathcal{O}_{S})\geq 0$. (Here we use notations in Theorem \ref{MT6} (ii.3).)
Since $\kappa(S)\geq 0$, we have $\chi(\mathcal{O}_{S})\geq 0$, that is, $A_{2}(X,L)\geq g_{1}(X,L)-1$.
Therefore we have $g_{1}(X,L)\leq 2$ from (\ref{MT5+1}).
If $g_{1}(X,L)=0$, then we have $h^{1}(\mathcal{O}_{X})=0$ by \cite[(4.8) Corollary and (1.1) Theorem]{Fujita89-3}.
But this contradicts the assumption that $h^{1}(\mathcal{O}_{X})>0$.
Next we assume that $g_{1}(X,L)=1$.
Then since $h^{1}(\mathcal{O}_{X})>0$, by \cite[(4.9) Corollary]{Fujita89-3} or Proposition \ref{T-P2} (iii) we see that
there exists a quasi-polarized variety $(V,H)$ such that $(V,H)$ is a scroll 
over a smooth elliptic curve and $(V,H)$ is birationally equivalent to $(X,L)$.
But this is impossible because we assume that $\kappa(K_{X}+2L)\geq 0$.
Therefore we have $g_{1}(X,L)=2$. In this case we see that $\chi(\mathcal{O}_{S})=0$.
\\
Since $\chi(\mathcal{O}_{S})=0$, we have $\kappa(S)=0$ or $1$.
\begin{Claim}\label{CL2}
$h^{1}(\mathcal{O}_{X})=2$.
\end{Claim}
\noindent
{\em Proof.}
Since $\dim\alpha(X)=2$, it suffices to show that $h^{1}(\mathcal{O}_{X})\leq 2$.
We also note that $h^{1}(\mathcal{O}_{X})=h^{1}(\mathcal{O}_{S})$.
\par
If $\kappa(S)=0$, then by the classification theory of surfaces we have $h^{1}(\mathcal{O}_{S})\leq 2$.
\par
So we may assume that $\kappa(S)=1$.
Then there exist a smooth curve $B$ and an elliptic fibration $\sigma:S\to B$.
In this case $h^{1}(\mathcal{O}_{S})\leq h^{1}(\mathcal{O}_{B})+1$ holds.
We consider the map $\psi:=\sigma\circ f:X^{\prime}\to S\to B$.
By taking the Stein factorization, if necessary, we may assume that $\psi$ is a surjective
morphism with connected fibers.
Let $F_{\psi}$ be a general fiber of $\psi$.
Since $\kappa(K_{X^{\prime}}+2\mu^{*}(L))=\kappa(K_{X}+2L)\geq 0$ by assumption,
we have $\kappa(K_{F_{\psi}}+2\mu^{*}(L)_{F_{\psi}})\geq 0$.
Hence by Proposition \ref{P1}, we get $h^{0}(K_{F_{\psi}}+2\mu^{*}(L)_{F_{\psi}})>0$.
Therefore $\psi_{*}(K_{X^{\prime}/B}+2\mu^{*}(L))\neq 0$ and
\begin{eqnarray}
h^{0}(K_{X^{\prime}}+2\mu^{*}(L))
&=&h^{0}(\psi_{*}(K_{X^{\prime}}+2\mu^{*}(L))) \label{MT5-2}\\
&\geq&\deg\psi_{*}(K_{X^{\prime}/B}+2\mu^{*}(L))
+h^{0}(K_{F_{\psi}}+2\mu^{*}(L)_{F_{\psi}})(g(B)-1). \nonumber
\end{eqnarray}
If $g(B)=0$, then $h^{1}(\mathcal{O}_{S})\leq 1$.
So we assume that $g(B)\geq 1$.
Since $1=h^{0}(K_{X}+2L)=h^{0}(K_{X^{\prime}}+2\mu^{*}(L))$,
we get $g(B)=1$ by Lemma \ref{L1}. Hence $h^{1}(\mathcal{O}_{S})\leq 2$.
\par
This completes the proof of Claim \ref{CL2}. $\Box$
\\
\par
By this claim, we see that $\alpha:X\to \mbox{Alb}(X)$ is surjective.
By \cite[Lemma 10.1 and Corollary 10.6]{Ueno75} we have $h^{2}(\mathcal{O}_{S})>0$ and $\kappa(S)=0$.
We also note that $\chi(\mathcal{O}_{S})=0$.
Hence $S$ is birationally equivalent to an Abelian surface.
Let $\tau: S\to S^{\prime}$ be the minimalization of $S$.
Then $S^{\prime}$ is an Abelian surface.
Here we note that there exists a rational map $\tau\circ f\circ\mu^{-1}:X\to S^{\prime}$.
By \cite[Lemma 9.11]{Ueno75}, this map is a morphism. We set $f^{\prime}:=\tau\circ f\circ\mu^{-1}$.
Let $F^{\prime}$ be a general fiber of $f^{\prime}$.
Then $F^{\prime}\cong \mathbb{P}^{1}$.
If $h^{0}(K_{F^{\prime}}+L_{F^{\prime}})>0$, then $h^{0}(K_{X}+L)>0$ by \cite[Lemma 4.1]{ChHa02}.
But this contradicts (\ref{MT5+2}).
Hence $h^{0}(K_{F^{\prime}}+L_{F^{\prime}})=0$.
Since $h^{0}(K_{F^{\prime}}+L_{F^{\prime}})=0$ and $h^{0}(K_{F^{\prime}}+2L_{F^{\prime}})>0$, we have $(F^{\prime},L_{F^{\prime}})=(\mathbb{P}^{1},\mathcal{O}_{\mathbb{P}^{1}}(1))$.
\\
(a.2.2) The case of $\alpha(X)=1$.
Then $\alpha(X)$ is a smooth curve and $\alpha: X\to \alpha(X)$ is a surjective morphism with connected fibers.
Let $C:=\alpha(X)$.
Here we note that $h^{1}(\alpha_{*}(K_{X}+2L))\leq h^{1}(K_{X}+2L)=0$
by the Leray spectral sequence.
Moreover $\kappa(K_{X}+2L)\geq 0$ implies $\kappa(K_{F_{\alpha}}+2L_{F_{\alpha}})\geq 0$
for a general fiber $F_{\alpha}$ of $\alpha$.
By Proposition \ref{P1} we have $h^{0}(K_{F_{\alpha}}+2L_{F_{\alpha}})>0$.
So we get $\alpha_{*}(K_{X/C}+2L)\neq 0$. 
On the other hand, by Lemma \ref{L1} we see that $\alpha_{*}(K_{X/C}+2L)$ is ample.
Hence $\deg \alpha_{*}(K_{X/C}+2L)>0$.
Therefore we have
\begin{eqnarray*}
h^{0}(K_{X}+2L)
&=&h^{0}(\alpha_{*}(K_{X}+2L)) \\
&=&h^{1}(\alpha_{*}(K_{X}+2L))+\deg\alpha_{*}(K_{X/C}+2L)
+h^{0}(K_{F}+2L_{F})(g(C)-1)\\ 
&\geq& 1+g(C)-1=g(C)\geq 1. 
\end{eqnarray*}

Since $h^{0}(K_{X}+2L)=1$, we have $g(C)=1$, that is, $h^{1}(\mathcal{O}_{X})=1$.
Hence $A_{2}(X,L)=g_{2}(X,L)+g_{1}(X,L)-h^{1}(\mathcal{O}_{X})=g_{2}(X,L)+g_{1}(X,L)-1$.
Since $g_{2}(X,L)\geq 0$ by Proposition \ref{2-P1}, we have $g_{1}(X,L)\leq 2$ from (\ref{MT5+1}).
By the same argument as (a.2.1) above, we get $g_{1}(X,L)=2$.
\par
We note that $h^{0}(K_{F_{\alpha}}+L_{F_{\alpha}})=0$ for a general fiber $F_{\alpha}$ of $\alpha$ because if $h^{0}(K_{F_{\alpha}}+L_{F_{\alpha}})>0$, then by \cite[Lemma 4.1]{ChHa02} we have $h^{0}(K_{X}+L)>0$ and this contradicts (\ref{MT5+2}).
Therefore by Proposition \ref{P1}  we have $\kappa(K_{F_{\alpha}}+L_{F_{\alpha}})=-\infty$.
Since $\dim F_{\alpha}=2$ and $\kappa(F_{\alpha})=-\infty$, we have $h^{0}(K_{F_{\alpha}}+L_{F_{\alpha}})=g(F_{\alpha},L_{F_{\alpha}})-h^{1}(\mathcal{O}_{F_{\alpha}})$ by the Riemann-Roch theorem and the Kawamata-Viehweg vanishing theorem.
Hence $h^{0}(K_{F_{\alpha}}+L_{F_{\alpha}})=0$ implies that $g(F_{\alpha},L_{F_{\alpha}})=h^{1}(\mathcal{O}_{F_{\alpha}})$.
By \cite[Theorem 3.1]{Fukuma97-1}, the $L_{F_{\alpha}}$-minimalization of $(F_{\alpha},L_{F_{\alpha}})$ (for the definition of the $L_{F_{\alpha}}$-minimalization see \cite[Definition 1.9 (2)]{Fukuma97-1}) is either $(\mathbb{P}^{2},\mathcal{O}_{\mathbb{P}^{2}}(1))$, $(\mathbb{P}^{2},\mathcal{O}_{\mathbb{P}^{2}}(2))$ or a scroll over a smooth curve.
But if $(F_{\alpha},L_{F_{\alpha}})\cong (\mathbb{P}^{2},\mathcal{O}_{\mathbb{P}^{2}}(1))$, then $h^{0}(K_{F_{\alpha}}+2L_{F_{\alpha}})=0$ and this is a contradiction.
Therefore we get the assertion of (a).
\\
(b) Assume that $h^{0}(K_{X}+3L)=1$.
Then by the proof of (a) in Theorem \ref{MT4-II}, we see that $\kappa(X)=-\infty$ and $A_{1}(X,L)\geq 1$.
But since $h^{0}(K_{X}+3L)=A_{1}(X,L)+2A_{2}(X,L)+A_{3}(X,L)$, $A_{2}(X,L)\geq 0$ and $A_{3}(X,L)\geq 0$, we see that $A_{1}(X,L)=1$.
Therefore $(g_{1}(X,L), L^{3})=(1,1)$ or $(0,2)$.
If $(g_{1}(X,L), L^{3})=(1,1)$ (resp. $(0,2)$), then we see that $(X,L)$ is the type (b.1) (resp. (b.2)) above by \cite[Corollaries (4.8) and (4.9)]{Fujita89-3}.\\
(c) Assume that $h^{0}(K_{X}+tL)={t-1\choose 3}$ for some $t\geq 4$.
Then by the proof of (b) in Theorem \ref{MT4-II}, we see that $A_{0}(X,L)=1$ and $A_{1}(X,L)=0$.
Hence $g_{1}(X,L)=0$ and $L^{3}=1$.
So we get the assertion by the same argument as in the proof of (a) in Theorem \ref{MT4-II}.
\par
These complete the proof. $\Box$

\begin{flushright}
Yoshiaki Fukuma \\
Department of Mathematics \\
Faculty of Science \\
Kochi University \\
Akebono-cho, Kochi 780-8520 \\
Japan \\
E-mail: fukuma@kochi-u.ac.jp
\end{flushright}


\begin{thebibliography}{99}

\bibitem{Andreatta95}
M. Andreatta,
\newblock{\em Some remarks on the study of good contractions,}
\newblock Manuscripta Math.\ 87\ (1995),\ 359--367.

\bibitem{AnWi97}
M. Andreatta and J. A. Wi\'sniewski,
\newblock{\em A view on contractions of higher dimensional varieties,}
\newblock Proc. Sympos. Pure Math. \ 62,\ Part 1\ (1997),\ 153--183.

\bibitem{BeSo95}
M. C. Beltrametti and A. J. Sommese,
\newblock{\em The adjunction theory of complex projective varieties,}
\newblock de Gruyter Expositions in Math.\ 16,\ Walter de Gruyter,\ Berlin,\ NewYork,\ (1995).

\bibitem{BlGe71}
S. Bloch and D. Gieseker,
\newblock{\em The positivity of the Chern classes of an ample vector bundle,}
\newblock Invent. Math.\ 12\ (1971),\ 112--117.

\bibitem{ChHa02}
J. A. Chen and C. D. Hacon,
\newblock{\em Linear series of irregular varieties,}
\newblock Algebraic Geometry in East Asia (Kyoto 2001),\ 143--153,\ World Sci. Publishing,\ River Edge,\ NJ,\ 2002.

\bibitem{Debarre}
O. Debarre,
\newblock{\em Higher-Dimensional Algebraic Geometry,}
\newblock Universitext, Springer-Verlag, New York, 2001.

\bibitem{Demaily96}
J. P. Demaily,
\newblock{\em Effective bounds for very ample line bundles,}
\newblock Invent. Math.\ 124\ (1996),\ 243--261.

\bibitem{EGA}
J. Dieudonn\'e and A. Grothendieck,
\newblock {\em \'El\'ements de G\'eom\'etrie Alg\'ebrique,}
\newblock Publ. Math. I.H.E.S.\ 4,\ 8,\ 11,\ 17,\ 20,\ 24,\ 28,\ 32.

\bibitem{E-V90}
H. Esnault and E. Viehweg,
\newblock{\em Effective bounds for semipositive sheaves and for the height of points on curves over complex function fields,}
\newblock Composit. Math.\ 76\ (1990),\ 69--85.


\bibitem{Fujita89-3}
T. Fujita,
\newblock{\em Remarks on quasi-polarized varieties,}
\newblock Nagoya Math. J.\ 115\ (1989),\ 105--123.

\bibitem{Fujita90}
T. Fujita,
\newblock{\em Classification Theories of Polarized Varieties,}
\newblock London Math. Soc. Lecture Note Ser.\ 155, Cambridge University Press,\ (1990).

\bibitem{Fukuma97-1}
Y. Fukuma,
\newblock{\em A lower bound for the sectional genus 
of quasi-polarized surfaces,}
\newblock Geom. Dedicata\ 64\ (1997),\ 229--251.

\bibitem{Fukuma98-2}
Y. Fukuma,
\newblock{\em On the nonemptiness of the linear system of polarized manifolds,}
\newblock Canad. Math. Bull.\ 41\ (1998),\ 267--278.

\bibitem{Fukuma04-1}
Y. Fukuma,
\newblock{\em On the sectional geometric genus of quasi-polarized varieties, I,}
\newblock Comm. Algebra\ 32\ (2004),\ 1069--1100.

\bibitem{Fukuma04-2}
Y. Fukuma,
\newblock{\em On the sectional geometric genus of quasi-polarized varieties, 
II,}
\newblock Manuscripta Math.\ 113\ (2004),\ 211--237.

\bibitem{Fukuma05}
Y. Fukuma,
\newblock{\em A lower bound for the second sectional geometric genus of polarized manifolds,}
\newblock Adv. Geom.\ 5\ (2005),\ 431--454.

\bibitem{Fukuma05-2}
Y. Fukuma,
\newblock{\em On the second sectional H-arithmetic genus of polarized manifolds,}
\newblock Math. Z.\ 250\ (2005)\ 573--597.

\bibitem{Fukuma06}
Y. Fukuma,
\newblock{\em On a conjecture of Beltrametti-Sommese for polarized $3$-folds,}
\newblock Internat. J. Math.\ 17\ (2006),\ 761--789. 

\bibitem{Fukuma07-2}
Y. Fukuma,
\newblock{\em On the dimension of global sections of adjoint bundles for polarized $3$-folds and $4$-folds,}
\newblock J. Pure Appl. Algebra\ 211\ (2007),\ 609--621.

\bibitem{Fukuma08-2}
Y. Fukuma,
\newblock{\em A study on the dimension of global sections of adjoint bundles for polarized manifolds,}
\newblock J. Algebra.\ 320 (2008),\ 3543--3558.

\bibitem{Fukuma10-2}
Y. Fukuma,
\newblock{\em On quasi-polarized manifolds whose sectional genus is equal to the irregularity,}
\newblock preprint (2010). http://www.math.kochi-u.ac.jp/fukuma/preprint.html

\bibitem{Horing09}
A. H\"oring,
\newblock{\em On a conjecture of Beltrametti and Sommese,}
\newblock preprint (2009), arXiv:0912.1295.

\bibitem{Horing09-2}
A. H\"oring,
\newblock{\em The sectional genus of quasi-polarised varieties,}
\newblock to appear in Arch. Math.

\bibitem{KMM85}
Y. Kawamata, K. Matsuda, and K. Matsuki,
\newblock{\em Introduction to the minimal model problem,}
\newblock Advanced Studies in Pure Math.\ 10\ (1985),\ 283--360.

\bibitem{Kleiman66}
S. L. Kleiman,
\newblock{\em Toward a numerical theory of ampleness,}
\newblock Ann. of Math.\ 84\ (1966)\ 293--344.

\bibitem{Mella}
M. Mella,
\newblock{\em Adjunction theory on terminal varieties,}
\newblock Complex analysis and geometry (Trento, 1995),\ 153--164,\ Pitman Res. Notes Math. Ser.,\ 366,\ Longman,\ Harlow,\ 1997.

\bibitem{Miyaoka87}
Y. Miyaoka,
\newblock{\em The Chern classes and Kodaira dimension of a minimal variety,}
\newblock Advanced Study in Pure Math.\ 10\ (1987),\ 449--476.

\bibitem{Ueno75}
K. Ueno,
\newblock{\em Classification Theorey of Algebraic Varieties and Compact Complex Spaces,}
\newblock Lecture Notes in Math.\ 439\ (1975), Springer.

\end{thebibliography}
\end{document}